\documentclass[12pt]{amsart}
\usepackage{amsbsy}
\usepackage{graphicx,epsfig,subfigure,psfrag}
\begin{document}

\newtheorem{theorem}{Theorem}
\newtheorem{proposition}{Proposition}
\newtheorem{lemma}{Lemma}
\newtheorem{corollary}{Corollary}
\newtheorem{definition}{Definition}
\newtheorem{remark}{Remark}
\newcommand{\tex}{\textstyle}
\numberwithin{equation}{section} \numberwithin{theorem}{section}
\numberwithin{proposition}{section} \numberwithin{lemma}{section}
\numberwithin{corollary}{section}
\numberwithin{definition}{section} \numberwithin{remark}{section}
\newcommand{\ren}{\mathbb{R}^N}
\newcommand{\re}{\mathbb{R}}
\newcommand{\n}{\nabla}
\newcommand{\iy}{\infty}
\newcommand{\pa}{\partial}
\newcommand{\fp}{\noindent}
\newcommand{\ms}{\medskip\vskip-.1cm}
\newcommand{\mpb}{\medskip}
\newcommand{\AAA}{{\bf A}}
\newcommand{\BB}{{\bf B}}
\newcommand{\CC}{{\bf C}}
\newcommand{\DD}{{\bf D}}
\newcommand{\EE}{{\bf E}}
\newcommand{\FF}{{\bf F}}
\newcommand{\GG}{{\bf G}}
\newcommand{\oo}{{\mathbf \omega}}
\newcommand{\Am}{{\bf A}_{2m}}
\newcommand{\CCC}{{\mathbf  C}}
\newcommand{\II}{{\mathrm{Im}}\,}
\newcommand{\RR}{{\mathrm{Re}}\,}
\newcommand{\eee}{{\mathrm  e}}
\newcommand{\LL}{L^2_\rho(\ren)}
\newcommand{\LLL}{L^2_{\rho^*}(\ren)}
\renewcommand{\a}{\alpha}
\renewcommand{\b}{\beta}
\newcommand{\g}{\gamma}
\newcommand{\G}{\Gamma}
\renewcommand{\d}{\delta}
\newcommand{\D}{\Delta}
\newcommand{\e}{\varepsilon}
\newcommand{\var}{\varphi}
\newcommand{\lll}{\l}
\renewcommand{\l}{\lambda}
\renewcommand{\o}{\omega}
\renewcommand{\O}{\Omega}
\newcommand{\s}{\sigma}
\renewcommand{\t}{\tau}
\renewcommand{\th}{\theta}
\newcommand{\z}{\zeta}
\newcommand{\wx}{\widetilde x}
\newcommand{\wt}{\widetilde t}
\newcommand{\noi}{\noindent}
\newcommand{\uu}{{\bf u}}
\newcommand{\xx}{{\bf x}}
\newcommand{\yy}{{\bf y}}
\newcommand{\zz}{{\bf z}}
\newcommand{\aaa}{{\bf a}}
\newcommand{\cc}{{\bf c}}
\newcommand{\jj}{{\bf j}}
\newcommand{\ggg}{{\bf g}}
\newcommand{\UU}{{\bf U}}
\newcommand{\YY}{{\bf Y}}
\newcommand{\HH}{{\bf H}}
\newcommand{\GGG}{{\bf G}}
\newcommand{\VV}{{\bf V}}
\newcommand{\ww}{{\bf w}}
\newcommand{\vv}{{\bf v}}
\newcommand{\hh}{{\bf h}}
\newcommand{\di}{{\rm div}\,}
\newcommand{\ii}{{\rm i}\,}
\newcommand{\inA}{\quad \mbox{in} \quad \ren \times \re_+}
\newcommand{\inB}{\quad \mbox{in} \quad}
\newcommand{\inC}{\quad \mbox{in} \quad \re \times \re_+}
\newcommand{\inD}{\quad \mbox{in} \quad \re}
\newcommand{\forA}{\quad \mbox{for} \quad}
\newcommand{\whereA}{,\quad \mbox{where} \quad}
\newcommand{\asA}{\quad \mbox{as} \quad}
\newcommand{\andA}{\quad \mbox{and} \quad}
\newcommand{\withA}{,\quad \mbox{with} \quad}
\newcommand{\orA}{,\quad \mbox{or} \quad}
\newcommand{\atA}{\quad \mbox{at} \quad}
\newcommand{\onA}{\quad \mbox{on} \quad}
\newcommand{\ef}{\eqref}
\newcommand{\ssk}{\smallskip}
\newcommand{\LongA}{\quad \Longrightarrow \quad}
\def\com#1{\fbox{\parbox{6in}{\texttt{#1}}}}
\def\N{{\mathbb N}}
\def\A{{\cal A}}
\newcommand{\de}{\,d}
\newcommand{\eps}{\varepsilon}
\newcommand{\be}{\begin{equation}}
\newcommand{\ee}{\end{equation}}
\newcommand{\spt}{{\mbox spt}}
\newcommand{\ind}{{\mbox ind}}
\newcommand{\supp}{{\mbox supp}}
\newcommand{\dip}{\displaystyle}
\newcommand{\prt}{\partial}
\renewcommand{\theequation}{\thesection.\arabic{equation}}
\renewcommand{\baselinestretch}{1.1}
\newcommand{\Dm}{(-\D)^m}

\title
{\bf  On regularity of a boundary point for higher-order parabolic
equations: towards Petrovskii-type criterion by  blow-up approach}


\author {V.A.~Galaktionov}

\address{Department of Mathematical Sciences, University of Bath,
 Bath BA2 7AY, UK}
\email{vag@maths.bath.ac.uk}



\keywords{Higher-order parabolic equations, boundary regularity,
 Petrovskii criterion,
 non self-adjoint operators, boundary layer, matching.
 {\bf Submitted to:}  NoDEA}

\date{\today}

\begin{abstract}

The classic problem of regularity of boundary points for
higher-order PDEs is concerned. For second-order elliptic and
parabolic equations, this study was completed by Wiener's (1924)
and Petrovskii's (1934) criteria,
 and
 was extended to more general equations including
 quasilinear ones.
Since the 1960--70s, the main success was achieved
for $2m$th-order elliptic PDEs (e.g.,  by Kondrat'ev and
  Maz'ya),
while the higher-order parabolic ones, with infinitely oscillatory
kernels, were not studied in such great details.


As a basic model, explaining typical difficulties of regularity
issues,
 the 1D bi-harmonic equation in a
 domain shrinking to the origin $(0,0)$ is concentrated upon:
  $$
  u_t=-u_{xxxx} \inB Q_0=\{|x| < R(t), \,\,-1<t<0\},
   $$
 where $R(t)>0$ is a smooth function on $[-1,0)$ and $R(t) \to 0^+$
as $t \to 0^-$. The zero Dirichlet conditions on the lateral
boundary of $Q_0$ and bounded initial data are posed:
 $$
 u=u_x=0 \atA x= \pm R(t), \,\, -1 \le t<0, \quad \mbox{and}
 \quad u(x,-1)=u_0(x).
  $$
The boundary point $(0,0)$ is then {\em regular} (in Wiener's
sense) if $u(0,0^-)=0$ for any data $u_0$, and is {\em irregular}
otherwise. The proposed asymptotic blow-up approach
shows that:

\ssk

  \noi (i) for the {\em backward  fundamental parabolae} with
$R(t)=l(-t)^{1/4}$, the regularity of its vertex $(0,0)$
depends on
 the constant
$l>0$:
e.g.,
 $l=4$ is  regular, while $l=5$ is
not;

\ssk

 \noi
 (ii) 
  for $R(t)=(-t)^{1/4} \var(-\ln (-t))$ with
$\var(\t) \to +\iy$ as $\t \to +\iy$,
 regularity/irregularity of $(0,0)$ can be expressed in
terms of
an integral Petrovskii-like (Osgood--Dini) criterion.
E.g., after a special ``oscillatory cut-off" of the boundary, the
function
 $$
 \tilde R(t) = 3^{-\frac 34}\, 2^{\frac {11}4}(-t)^{\frac 14} \big[\ln |\ln(-t)|\big]^{\frac 34}
 $$
 belongs to the regular  case, while any increase of the constant
 $ 3^{-\frac 34}\, 2^{\frac {11}4}$ therein
leads to the irregular one.
 The results are based on Hermitian
spectral theory of the operator
 $
 \BB^*=- D_y^{(4)}- \frac 14 \,\, y D_y$ in
 $L^2_{\rho^*}(\re)$, where
 $\rho^*(y) = {\mathrm e}^{-a |y|^{4/3}}$, $a={\rm constant} \in (0,3 \cdot
 2^{-\frac {8}3})$,
together with typical ideas of boundary layers and blow-up
matching analysis.
Extensions to $2m$th-order poly-harmonic equations in $\ren$ and
other
  PDEs are
 discussed,
  and a partial survey on regulariy/irregularity issues
 is presented.

\end{abstract}

\maketitle

\begin{center}
{\em Dedicated to  the memory of Professor
I.G.~Petrovskii}
\end{center}

\section{Introduction: Petrovskii's regularity  criterion
 (1934) and  extensions}
\label{S1}

\subsection{First discussion: regularity as a fundamental issue in potential theory}

Well-posedness of
 initial-boundary value problems (IBVPs) for linear and
nonlinear partial differential equations (PDEs) are key for
general PDE theory and applications. The principal question, which
was always in the  focus of the key research
  in this  area from
the nineteenth century to
 our days,
 is to
determine
 optimal, and as sharp as possible, conditions on the ``shape"
 of the continuous boundary, under which the solution is continuous at
 the boundary points. This means that, for a given PDE with
 sufficiently smooth coefficients, a standard IBVP can be posed with a classical notion of solutions
  applied.

 It is impossible to mention all the cornerstones of such
 regularity PDE theory achieved in the twentieth century.
In the nineteenth century, this study, for the {\em Laplace
equation}
 \be
 \label{Lap1}
 \D u=f \inB \O \subset \re^2 \,\,\, \mbox{or}\,\,\, \re^3
 \quad (0 \in \partial \O \,\,\,\mbox{is the point under scrutiny}),
  \ee
 began by Green (since 1828), Gauss (1940), Lord Kelvin
(1847) and Dirichlet himself (in the 1850s), Weierstrass and
Neumann (in the 1870s), Hilbert (1899), Schwarz, Poincar\'e (since
1887); see detailed history of potential theory in Kellogg
\cite[pp.~277--285]{Kell29}.
It is also worth mentioning that
classic in potential theory ``Lyapunov's surfaces" were introduced
   about  1898
 for general non-convex domains in $\re^3$ (in the
 convex case, the problem was  solved  in 1870 by Karl
 Neumann). Lyapunov's study was known to be inspired by Poincar\'e's earlier
 papers;
 a full proof for general Lyapunov surfaces was
 completed by V.A.~Steklov about 1902. As a related issue, as was pointed
 out by V.G.~Maz'ya \cite{Maz08},
 Poincar\'e seems was the first who in last years of the nineteenth  century already knew and
  used the exterior
 cone condition for the regularity of boundary points. As is well known from his
 private letters,  discussions with several Russian mathematicians, including
 e.g., his supervisor P.L.~Chebyshov\footnote{A.M.~Lyapunov Master's Thesis ``On Stability of
Spheroidal Equilibrium Forms of a Rotating Fluid", written under
the supervision of Chebyshov, was defended in the S.-Petersburg
University, 27th January, 1885.}
  and Steklov, Lyapunov was constantly rather anxious studying
new and amazing but not fully   rigorously justified  ideas and
methods by Poincar\'e. Possibly, this led him to develop the
concept of ``Lyapunov's surfaces" in potential theory about that
time (as well as some of other fundamental concepts in his
stability theory, where discussions with Poincar\'e were also
known to take place).

 For the second-order Laplace equation \ef{Lap1} and the  {\em heat equation}
 \be
 \label{HE11}
  u_t= \D u \inB Q, \quad (0,0) \in \partial Q \,\, \mbox{is
  characteristic, where $\partial Q$ is tangent to $\{t=0\}$},
   \ee
 regularity theory  was almost fully
  completed in the 1920--30s, and
  various ideas  and key results
  on
  optimal regularity criteria
 created  the amazing history,
 which is explained in a number of classic monographs.
For convenience,  we will begin our discussion with some
historical aspects  in next Section \ref{S1.2}, where we hope to
present some new and not that well-known regularity features and
peculiarities for the sake of the attentive Reader.

 The present paper
 uses the asymptotic blow-up approach from reaction-diffusion
 theory
 to such optimal boundary regularity questions.
  It is shown that, for a number of higher-order PDEs with essentially oscillatory
kernels,
  these can be
  treated by  {\em blow-up evolution} via approaching a
 ``singular" boundary point.
Surely, this is not a principal novelty not only for the heat
equation \ef{HE11} (Petrovskii, 1934), but also in  elliptic
regularity theory, where, together with other techniques, rescaled
linear operator pencils were shown to be key; we refer to
Kozlov--Maz'ya--Rossmann's monographs
 \cite{KMR1, KozMaz01} and surveys in \cite{Freh06, KondOl83, Maz02,  Wood07}
  for details and full history.
In this connection, Kondrat'ev's seminal paper \cite{Kond67} for
higher-order linear elliptic equations\footnote{As explained in
\cite{Kond67}, previous results on normal solvability of elliptic
boundary value problems in domains with angular (in $\re^2$) and
conical ($\ren$, $N \ge 3$) points were obtained by Fufaev ($m=1$,
$N=2$, 1960), Volkov (1962), Birman and Skvortsov (1962),
Lopatinski$\breve{\rm i}$ (1963). Earlier Nikol'skii's results
\cite{Nik56} in the 1956--58 on boundary differential properties
of functions in the Nikol'skii spaces $H^r_p$ defined on regions
in $\re^2$ with angular points and hence criteria of such a
regularity for the Laplace equation should be mentioned.}
 (see  Kondrat'ev--Oleinik's very detailed survey of 1983 \cite{KondOl83}
for up to the 1980s and the most recent monograph \cite{Bor07} for
a huge number of more modern extensions), thought not devoted to
the regularity issues, represented a novel and involved use of
spectral properties of operator pencils for describing the
asymptotics of solutions near corner singularities. This analysis
assumes a kind of  an ``elliptic evolution" approach for elliptic
problems, which then is not well-posed in Hadamard's sense but can
trace out the behaviour  of necessary global orbits that approach
the singularity point. For nonlinear elliptic equations, these
ideas got various developments and diversions since the 1970s; see
V\'eron's monographs
 \cite{VerMon, VerMon2} and \cite{BidVer99} for further history and key
references.
 Classic Wiener's and operator pencil ideas were successfully developed for many higher-order elliptic
 equations with a typical representative:
  \be
  \label{Rep1}
(-\D)^m u=f \inB \O \subset \ren \quad (m \ge 2),
 \ee
in the papers by Maz'ya with collaborators since the 1970s; see
e.g., \cite{MazPl81} and \cite{Maz99, Maz02, KMR1, KozMaz01} as a
guide. Though regularity point analysis for \ef{Rep1} revealed a
lot of new phenomena and mathematical difficulties relative the
classic ones for $m=1$ (see comments later
 on), there are many papers on this subject, where essential
 progress was achieved, especially in the recent ten years.

\subsection{Higher-order parabolic PDEs: on some known regularity results}


For {\em $2m$th-order parabolic equations}, where the basic model
is the canonical {\em poly-harmonic equation}
 \be
 \label{PE1}
 u_t=-(-\D)^m u \withA m \ge 2,
  \ee
regularity questions have been also addressed in the literature,
however essentially less actively  than for elliptic PDEs referred
to above. Among results of classic solvability theory, which can
be found in a number of monographs on linear and nonlinear
equations, systematic approaches to the study of non-cylindrical
domains with a characteristic vertex $(0,0)$,
 where\footnote{Here and later on, we indicate only the results
 that are related to our blow-up setting to be explained shortly;
 these deep papers contain other involved conclusions.}
  the tangent
plane $\{t={0}\}$ is horizontal in the $\{x,t\}$-space appeared
already in the 1960s. V.P.~Miha${\rm \check{i}}$lov papers in 1961
and 1963 \cite{Mih61, Mih63I, Mih63II} treated the case, where the
parabolic boundary of the domain $Q_0$ lies below the
characteristic plane $t=0$ and has a typical form of the
corresponding ``backward" paraboloid of the fundamental
solution\footnote{{\em Q.v.} more recent Miha${\rm
\check{i}}$lov's research on existence of boundary values of
poly-harmonic functions for domains with smooth boundary; see
references in \cite{Mikh05}.}
 (in \cite{Mih63I, Mih63II}, this shape is
perturbed by factor that makes the paraboloid ``sharper", i.e.,
becomes ``more regular"; see further comments on that).
 Almost simultaneously with the fundamental
``elliptic" paper \cite{Kond67}, Kondrat'ev in 1966 published
another key ``parabolic" paper \cite{Kond66}, which  is less
known\footnote{It is useful to compare the total number of
citations in the {\tt MathSciNet} (October 2008): the elliptic
paper \cite{Kond67} has the record 228 citations, while the
parabolic one \cite{Kond66} has 7 citations, i.e., 32 times less!
However, the proposed scaling techniques \cite{Kond66, Kond67}
leading to operator pencils are almost identical in both parabolic
and elliptic papers (moreover,  the parabolic pencil analysis
\cite{Kond66} directly applies to some ultraelliptic PDEs). The
asymptotic expansions derived by spectral theory near the vertex
0,
 \be
 \label{KondEx}
  \tex{
 u(x,t) \sim \sum_{(j,k)} \t^{-\ii \l_j/2m} (\ln \t)^k \,
 \psi_{kj}
  \quad (\l_j\,\, \mbox{are eigenvalues of
 pencils})
 }
  \ee
 are similar for the elliptic ($\t=|x|$) and parabolic ($\t=-t$)
 cases (the $\ln^k \t$, $k=0,1,...,k_j$ multipliers reflect finite algebraic multiplicities
  $k_j$ of $\l_j$
 with associated generalized $C^\iy$ eigenvectors $\{\psi_{kj}\}$).
 Note that, even for the heat equation, the required pencil spectral properties are not that easy to
 get, \cite{Ar98}.}. Here, Kondrat'ev
 dealt  with the singularity issues
 for more general than \ef{PE1} parabolic equations and various boundary
 conditions on a fixed  fundamental paraboloid. Complicated asymptotics of
 solutions near such singular points were obtained, and, as a
 result,
 a solvability criterion was proved for
existence of  weak solutions  in special weighted
Sobolev--Slobodetskii-type spaces. In 1971, Fe${\rm \check{i}}$gin
\cite{Fei71} extended the results to the case when, in any
neighbourhood of the vertex $(0,0)$, the boundary lies on both
sides of the characteristic plane $\{t=0\}$ (see \cite{Rab04} for
Fe${\rm \check{i}}$gin's contribution to elliptic regularity
theory in the 1970s). These deep results initiated further
subsequent study, which nevertheless was not that exhaustive and
sharp as in elliptic theory, and we will try to explain why (the
answer seems easy: another level of difficulty). The current
standing of this regularity parabolic theory can be traced out
further  by the {\tt MathSciNet}, but the author reports that
getting a clear convincing view  is not  easy.


\subsection{Main goal: deriving sharp asymptotics of solutions near
vertex depending on boundary geometry}

As a consequence of the above discussion involving a number of
strong previous results,
obviously,
  optimal conditions of regularity/irregularity  of the vertex $(0,0)$
   for the poly-harmonic flow \ef{PE1} demand the following:
 \be
 \label{G001}
  \begin{matrix}
 \mbox{{\bf Goal I}: \,\, sharp asymptotic expansion of solutions
 near $(0,0)$ are necessary}\qquad \ssk\\
 \mbox{for rather arbitrary backward characteristic (non-fundamental) paraboloids}
 \ssk\\
 \mbox{to get regularity ot the vertex {\em via} the
 geometry of the boundary nearby.}
 \qquad
  \end{matrix}
  \ee
 In other words, we check whether it is possible to
characterize a sharp asymptotic behaviour of classic (strong)
solutions near the characteristic vertex, when the paraboloid {\bf
is not} given by the backward variable of the fundamental
solution, i.e., can take various forms.
 Loosely speaking, we are going to  mimic   a Petrovskii-type
 criterion for higher-order parabolic equations such as \ef{PE1}
 (we must admit that, literally, this goal is not  and
 cannot be achieved literally, but a partial progress is indeed doable).
Once such an asymptotic has been obtained, one can re-evaluate the
solution near the vertex in any necessary weighted metric,
eventually to check whether solution exists (is continuous at the
vertex) in weaker sense.

As usual, we find useful, in order to characterize our blow-up
approach, to consider first
the simplest
 1D {\em
 bi-harmonic equation}
  \be
  \label{bih1}
  u_t=-u_{xxxx} \quad (m=2),
   \ee
 which even in the one $x$-dimension  provides us with a number of not that mathematically pleasant surprises.
 It seems that the main difference (and the origin of extra difficulties)
  between the higher-order elliptic case
 \ef{Rep1} and the parabolic one \ef{bih1} is that the later
  has always  oscillatory kernel of changing sign infinitely many
 times, while for \ef{Rep1}, the kernel is   less oscillatory
 and even remains positive in some dimensions even for $m \ge 2$.
 As is well known since Maz'ya--Nazarov's  results for $m=2$, $N \ge 8$ in
 1986
 \cite{MN86}, the changing sign kernels of elliptic operators
 produce new phenomena, where classic regularity techniques may
 fail and even standard conical points can be irregular (and singular, i.e., unbounded at the vertex).
 Note that Kondrat'ev's sharp estimates such as \ef{KondEx}
 \cite{Kond67} can also indicate such a possibility, if sharp
 bounds on first pencil's eigenvalues 
  are known.
Therefore, a simple 1D model \ef{bih1} becomes a key parabolic PDE
with such an infinitely oscillatory kernel, so we concentrate on
this equation in what follows.

 As a unified issue of the present analysis of higher-order parabolic and other PDEs (i.e., without any
 order-preserving or positive kernel
 features), we aim a sharp
  treatment of boundary regularity conditions
  by
  using  spectral theory of related rescaled non
 self-adjoint operators or pencils of such operators.
  A full
 justification of some ``approximate" regularity conditions are very
 difficult and lead to open problems.
 As an independent feature of the regularity questions, which
 underlines their  difficult nature,
  as far as we know, even  for the simplest higher-order elliptic or parabolic
 operators with real constant coefficients,
  {there is a clear and natural  difficulty for determining
     sufficiently sharp regularity conditions for
  linear operators with kernels of essentially changing sign.}
Of course, we exclude those operators for which the regularity is
guaranteed by standard embedding results, or when the irregularity
is induced by some symplectic--Hamiltonian properties, which do
not allow any shrinking of the domain under consideration (e.g.,
the latter is true for  Schr\"odinger operators\footnote{The
author thanks I.V.~Kamotski for a discussion that clarified this
 property.}
 $\ii D_t +
(-\D)^m$, where $L^2_x$-norm is preserved).
 For such
operators,
Wiener's and other related capacity-like techniques, which assume
kernel positivity, are naturally expected to fail, so new
techniques are necessary.

\ssk

Thus,  our {\em goal} is to develop an asymptotic method of
blow-up regularity analysis, firstly, for the simplest 1D {\em
bi-harmonic equation}, and next extend to poly-harmonic ones
\ef{PE1}. In general, together with other examples of PDEs,
we follow the principle (as we have seen, such kernels have been
studied much less in the mathematical literature)
\be
\label{OscInf}
 \mbox{{\bf Goal II}: \,\,to deal with kernels having infinitely many sign changes,}
  \ee
  that can occur not necessarily for parabolic PDEs only.
We will show that a kind of Petrovskii-like integral ``criterion"
for regularity/irregularity of boundary points \cite{Pet34Cr,
Pet35} (1934--35) can be visualized.
However, {for operators with oscillatory kernels of changing sign,
a standard deterministic
 analysis of such integral Osgood--Dini-type conditions via
 divergence/convergence
(as for the second-order or other operators with more positively
 dominant kernels) is not available.}
 Secondly, we discuss possible
extensions of the method to other types of linear and nonlinear
PDEs with the kernels like \ef{OscInf}. Meantime, we continue our
survey of the history of crucial for us regularity issues and {\em
integral criteria}.


 \section{Elliptic equations as the origin of regularity theory:
  Zaremba (1911),  Lebesgue (1913), and Wiener (1924)}
 \label{SEll1}




\subsection{Laplace equation: early days
of existence and nonexistence theory}

  The first {\em criterion} (a necessary and sufficient condition) of
regularity of a boundary point 0 was the
 Wiener famous one  for the Dirichlet problem for the {Laplace
 equation}
 $$
  \D u=f \,\,\, \mbox{in a bounded domain}\,\,\, \O \subset \re^2,
 \,\,\, 0 \in \partial \O \quad \big(u=0 \,\, \mbox{on $\partial \O$,
 for definiteness}\big),
 $$
derived in 1924, \cite{Wien24} (the case of  arbitrary $N \ge 2$
was also embraced).
 It is formulated in terms of a diverging
 series of capacities (measuring the thickness of the complement
 of $\O$ near 0) of a discrete family
 of domains shrinking to the given point $0$ on the boundary;
 see extra details and further discussions  in Courant--Hilbert \cite[p.~306]{CourH} and in Kellogg
 \cite[p.~330]{Kell29}.

 Concerning {\em nonexistence} for the Dirichlet problem
 (see \cite{Maz97, KMR1, KozMaz01} for details),
the first rigorous nonexistence example belonged to Stanislav
Zaremba  in 1911 \cite[p.~310]{Zar10} (for a plane domain whose
boundary has isolated points, not connected;
  e.g., the surface comprised a sphere
  and its centre); see comments in \cite[p.~285]{Kell29}.
 The key non-solvability example was constructed by
 Lebesgue in 1913
 \cite{Leb13, Leb24}\footnote{A full Lebesgue account on the Dirichlet
 problem includes his important earlier paper \cite{Leb07}, where
 the variational method was employed in 2D to solving under
 very general conditions, and another key one \cite{Leb12Bar},
 where ``barrier" (in Poincar\'e's terminology) techniques were employed; \cite{Leb24} is his final summarizing all key
  ideas and results.}
 who discovered striking  examples
   of irregularity. E.g., he showed that,
for a domain in $\re^3$ with the {\em inner cone}  obtained by
rotating the function $x_2={\mathrm e}^{-1/x_1}$ for $x_1>0$
small, about the $x_1$-axis (then the points with $x_2>{\mathrm
e}^{-1/x_1}$ do not belong to the domain, i.e., the domain has a
sharp cusp  turned inwards, a {\em Lebesgue spine}), the origin 0
is irregular. Figuratively speaking, despite such a thin spine
connection of the 0 with the boundary,  ``effectively" the centre
remains ``separated" reminding earlier Zaremba's nonexistence
construction.
 The same is true for $x_2=x_1^{|\ln x_1|^\e}$, with any
$\e>0$ (0 is regular if it can be touched from outside by any cone
given by rotation of $x_2=x_1^k$ for any $k>0$ fixed); see
Petrovskii's text book \cite[p.~325]{Pet61}.
 Carleman's PhD Thesis \cite{Carl16} of 1916 was also one of the
 first study of singularities of solutions of elliptic PDEs at non-regular boundary
 points; see  also \cite[Introd.]{KondOl83}.

 Forty years later,
 new
   examples of nonexistence (somehow in lines with  Zaremba--Lebesgue--Uryson's ideas)
   were developed in \cite[p.~129]{Kond61} by Kondrat'ev.
   For the Dirichlet problem $\D u+u=0$ in $\O$,
he produced further examples of non-solvability for not simply
connected domains for $N=2$ (sufficient number of disjoint small
discs are removed), and  domains in $\re^3$ such as a ball   with
line segments or narrow cylinders removed to get sufficiently
small interior diameter (one end stays on the surface for not
violating connectedness). These constructions extend to arbitrary
$N$ and $2m$, including the most delicate range $N>2m$.
 By an $m$-capacity technique,
\cite{Kond61} also
studied the question of $k$-regularity ($k \le m$) at a boundary
point $A$ meaning that
any $u \in
 W^2_{m,0}$ has $k-1$ derivatives vanishing there {\em asymptotically}
 (up
  to a set of the measure
  $o(a^N)$ in $|x-A|<a$ as $a \to 0$).





\subsection{Uryson (1924), Keldysh (1941), and Miusskaya
Square}

 A nonexistence example similar to Lebesgue's one  was
independently constructed by P.S.~Uryson (1898--1924) in the early
1920s; his paper \cite{Urys25} was published in 1925, i.e., after
his tragic death. From electrostatic point of view, which is
natural for such problems, the electricity ``flows down" from the
cusp, i.e., cannot be retained there.
Besides taking a physics curse from his PhD adviser N.N.~Luzin
(Luzin's postgraduate course was completed by him in 1921), the
physical intuition of Uryson was possibly created even before,
when he, being in a gymnasium in Moscow,  began to work at the
Shanyavskii University in the Miusskaya Square, Moscow, under the
supervision of the eminent Russian physicist P.P.~Lasarev, and
even published a paper on experimental research in X-radiation;
 see \cite{UrOb98} (it seems less known that Uryson was the first,  before Luzin, supervisor
 of A.N.~Kolmogorov in his student's period). This emphasizes Uryson's
  outstanding physical motivation and understanding
 in the very early  days of his scientific carrier.
  Both the Shanyavskii University and the Physical
 Institute, founded in 1911 for the outstanding physicist P.N.~Lebedev (who first
 measured ``pressure of the sun light"),
were erected in 1910--1912 by Russian architect I.I.~Ivanov-Schitz
 by using funding from Moscow's merchants\footnote{This old building (as we
 will see, somehow related to  boundary regularity) can be seen in
 ``The World of Andrey Sakharov" (he first visited it
 in 1945), in {\tt http://people.bu.edu/gorelik/Publications}.}.
Later on,  since 1953, on the basis of  Lebedev's Institute
building, was developed the Keldysh
Institute of Applied Mathematics\footnote{The host institution for
the author of this paper for almost twenty five years (as an
excuse for paying possibly too much unreasonable attention to
those historical events around it).}
 founded and headed until his death in 1978 by
  M.V.~Keldysh.

   Interestingly (and key for the whole comment),
   Keldysh made a significant and novel
contribution to Dirichlet's problem theory in 1941 in his
celebrated paper
\cite{Keld41} (see also \cite{Keld41CR}), which  treated irregular
points, as well as the new notion of stability of points and of
the problem inside the domain relative perturbations of the
boundary (the notion of capacity is discussed in Ch.~II, while
Ch.~III is devoted to Wiener's regularity criterion).  In Ch.~V,
he proved that a stable boundary point is always regular, but the
converse is not true
 ({\em Keldysh balls} $M$ were constructed having  zero area of unstable points on
$\partial M$ and of positive harmonic measure; see \cite{Ishi97,
Ar08} for a modern exhibition and extension).

 It is also worth
mentioning that ten years later,  in 1951, Keldysh published the
related paper \cite{Keld51Deg} on elliptic equations that are
degenerated on the boundary, where particularly  the novel idea
that a part of the boundary may be free from any conditions was
introduced the first. These ideas were later developed by Fichera
(since 1956), Oleinik, Radkevich, Kohn, Nirenberg (in the 1960s),
etc.; see the monograph \cite{OlRad73} and \cite{Chun07} for a
more recent standing. For a full collection, let us mention
Keldysh's fundamental paper \cite{Keld51Pen} (same 1951 year!) on
linear operator pencil theory {\em with  applications to elliptic
PDEs}. Therein, several novel ideas including difficult $n$-fold
completeness questions of root functions, i.e., eigen and
associated ones, for the {\em Keldysh pencils} with variable
coefficients, {\em Keldysh chains}, etc. were introduced. Its
sequel in the 1950s includes works by I.M.~Glazman, M.G.~Krein,
V.B.~Lidskii, and others; these Keldysh methods of
\cite{Keld51Pen} were later used by Browder, Agranovich,
Krukovsky,  Agmon, Egorov, Kondrat'ev, Schulze, etc.; see
\cite{Eg02} for details. The paper \cite{Keld51Pen}  was later
classified as one of the main achievements of pencil theory in the
twentieth century (together with J.D.~Tamarkin's results in his
PhD Thesis in 1917, S.-Petersburgh's University \cite{Tam17}); see
Markus' monograph \cite{Markus} for the history and key results.



\subsection{On some  recent achievements}


  Since we are not dealing
with elliptic equations in what follows, it is worth mentioning
here that a full extension of Wiener's regularity test (criterion)
via the concepts of potential-theoretic Bessel (Riesz) capacities
to $2m$th-order strongly elliptic equations with real constant
coefficients,
 \be
 \label{LL1}
  \left\{
  \begin{matrix}
 L(\partial)u=f, \quad f \in C_0^\iy(\O), \quad u \in  H^m_0,
 \qquad\qquad\qquad\qquad\qquad\qquad\quad  \ssk\ssk\\
L(\partial)= \sum\limits_{|\a|=|\b|=m}a_{\a\b} \partial^{\a+\b},
\,\, a_{\a\b} \equiv a_{\b\a}, \,\,\, (-1)^m L(\xi)>0 \,\,
\mbox{in}\,\, \ren \setminus\{0\},
 \end{matrix}
 \right.
  \ee
  was completed in the case $N=2m$
     in 2002
   by Maz'ya\footnote{The principal extension of Wiener's-like capacity
   regularity  test to the nonlinear degenerate $p$-Laplacian operator
   $-\D_pu \equiv -\n \cdot(|\n u|^{p-2} \n u)$, $p \in
   (1,N]$, was also due to Maz'ya in 1970 \cite{Maz70}; later on,
   his sufficient capacity regularity condition turned out to be optimal for any $p>1$.} ($N>2m$ also admits a treatment,
   though applied for a subclass of the so-called {\em strongly
   elliptic equations};
  for $N \le 2m-1$, the point is always regular by the classic
Sobolev embeddings); see \cite{Maz99, Maz02}
  and references therein to earlier results and extensions.



 Concerning the related questions of  sharp asymptotic behaviour near Lipschitz boundary
at 0 for problems \ef{LL1}, in \cite{KozM05},
  a rather general case of the boundary  given by
  $$
  \O=\big\{x=(x',x_n) \in \ren: \quad x_n> \varphi(x')\big\}, \quad
  \varphi(0)=0,
   $$
is studied.  The asymptotics as $x \to 0$ of a class of solutions
are shown to depend on the Lipschitz constant of $\varphi(y')$ in
a neighbourhood of 0:
 $$
  \tex{
  \kappa(\rho)= \sup_{|y'|<\rho}\, | \n \varphi(y')|
  \quad\big(\mbox{estimates are simplified if} \,\,\,
   \int_0 \frac{\kappa^2(\rho)}\rho\, {\mathrm d}\rho < \iy\big).
   }
   $$

    Among
  others, let us also note an unusual Wiener's-type regularity and existence result for {\em
  blow-up on the boundary} solutions  of the semilinear elliptic equation
  \be
  \label{DD1}
 \D u-|u|^{q-1}u=0 \inB \O, \quad u(x) \to + \iy \asA x \to
 \partial \O \quad (q>1);
   \ee
 see \cite{Lab03, Marc07}, where a long list of other related references
 is presented (equation \ef{DD1} will be cited again).



\section{Parabolic PDEs:
 Petrovskii's $\sqrt{\log \log}$ criterion (1934) and around}
 \label{S1.2}

 \subsection{Parabolic regularity theory began with the
heat equation: first results and definitions}

 From the beginning, for the {\em second-order parabolic
equations}, the boundary point regularity analysis took a
different direction (i.e., not of Wiener's type).
  This was Petrovskii \cite{Pet34Cr, Pet35}, who in
  1934\footnote{Actually, this research was performed earlier: his most famous 1935 paper in
  {\tt Compositio Mathematica} \cite{Pet35} was submitted for publication
 ``(Eingegangen den 27. November 1933.)" \cite[p.~419]{Pet35}. Before, the results were
 discussed,
  ``Diskussion im wahrscheinlichkeits-theoretischen Seminar der Universit\"at Moskau",
   which led to Kolmogorov's problem solved
  by Petrovskii in \cite[p.~414]{Pet35}; see Abdulla \cite{Abd05, Abd08Adv} for further details concerning this
  fundamental problem called nowadays  Kolmogorov--Petrovskii's one.}
   completed the study
  of
the regularity question for the 1D and 2D {\em heat equation} in a
non-cylindrical domain.
 For further use, we formulate his result in a
blow-up manner, which in fact was already  used by Petrovskii in
1934 \cite{Pet34Cr}. This is about the question on a {\em
irregular} or {\em regular} point $(x,t)=(0,0)$ for nonnegative
solutions of the IBVP
 \be
 \label{ir1}
  \left\{
 \begin{matrix}
   u_t=u_{xx} \quad \mbox{in}
 \quad Q_0=\{|x|<R(t), \,\,\, -1 < t<0\}, \quad R(t) \to 0^+ \,\,\mbox{as}\,\, t \to
 0^-, \ssk\ssk\ssk\\
  \mbox{with bounded smooth data $u(x,0)=u_0(x)\ge 0$ on $[-R(-1),R(-1)]$}.
  \quad\,\,\,
   \end{matrix}
    \right.
  \ee
Here the lateral boundary $\{x= \pm R(t), \,\, t \in [-1,0)\}$ is
given by a function $R(t)$ that is  assumed to be positive and
$C^1$-smooth for all $-1 \le t  <0$ and is allowed to have  a
singularity of $R'(t)$ at $t=0^-$ only. We then study the value of
$u(x,t)$ at the end ``blow-up" {\em characteristic}  point
$(0,0^-)$, to which the domain $Q_0$ is said to ``shrink" as $t
\to 0^-$.

 For the heat equation, the first existence  of a classical solution result was obtained  by
 Gevrey in 1913--14 \cite{Gev13}
 (see Petrovskii's references  in \cite[p.~55]{Pet34Cr} and \cite[p.~425]{Pet35}),
 which assumed that the H\"older exponent of
 $R(t)$ is larger than $\frac 12$, i.e., in our setting, at $t=0^-$, all types
 of boundaries given by the functions
  \be
  \label{Gev1}
  \fbox{$
  R(t)=
(-t)^\nu \quad \mbox{with any $\nu > \frac 12$ are regular}
 \quad(\mbox{Gevrey, 1913--14}).
 $}
 \ee
Note that, for $2m$th-order parabolic equations, a similar result
saying that, for
 \be
  \label{Gev1Mih}
  \fbox{$
  R(t)=
(-t)^{\frac 1{2m}}, \quad \mbox{the problem is uniquely solvable}
 $}
 \ee
(in Slobodetskii--Sobolev classes, i.e., not in the classic sense)
was proved by Miha$\breve{\rm i}$lov \cite{Mih61} almost sixty
years later.

\ssk


\noi\underline{\bf Regular point}: as usual in potential theory,
the point $(x,t)=(0,0)$ is called {\em regular} (in Wiener's
sense, see \cite{Maz99};
 sometimes is also called {\em
boundary}) if any value of the solution $u(x,t)$ can be prescribed
there by continuity as a standard boundary value on $\partial
Q_0$. In particular,
 as a convenient and key for us  evolution illustration, $(0,0)$ is regular if the
 continuity holds for any initial data $u_0(x)$ in the following sense:
  \be
  \label{u0}
  u=0 \,\,\, \mbox{at the lateral boundary} \,\,\,\{|x|=R(t), \, -1 \le t<0\} \LongA
  u(0,0^-)=0.
  \ee

\noi\underline{\bf Irregular point}: otherwise, the point is {\em
irregular} (or {\em inner}) if the value $u(0,0^-)$ is not
arbitrary and is given by the parabolic
evolution as $t \to 0^-$ (hence, formally and actually, $(0,0)$
does not belong to the parabolic boundary of $Q_0$, i.e., is {\em
inner} in a natural sense).





\subsection{On some details of Petrovskii's analysis in 1934--35 and extensions}

Reading last twenty five or so years various papers related to
regularity issues for parabolic PDEs, the author found several and
sometimes rather distinctive treatments of Petrovskii's results.
For instance, his 1934 paper \cite{Pet34Cr} was not cited
practically in no such papers, but the one \cite{Ar98}. There were
also discrepancies in treating whether or not Petrovskii derived
the integral Osgood--Dini-like conditions (surely, he did already
in 1934).
 Therefore, we expect that it is worth mentioning
here some aspects of Petrovskii's original derivation of his
celebrated criterion.

Thus, according to \ef{ir1},  in both cases of
regularity/irregularity, there occur interesting asymptotic
problems of the behaviour of $u(x,t)$ in $Q_0$ as $t \to 0^-$,
which are of singular (finite-time blow-up) type.
 These problems were solved by Petrovskii \cite{Pet34Cr, Pet35} in
 1934--35,
 who derived a regularity criterion
\footnote{Cf. quoting these necessary and sufficient conditions
  in Kondrat'ev
 \cite[p.~450]{Kond66}.}
   for \ef{ir1}
 by constructing rather tricky sub- and super-solutions
 (``subparabolish" and "superprabolish" in the original Petrovskii
 terminology, \cite[p.~386]{Pet35})
 of the heat equation.

  For instance,  to prove irregularity for the characteristic
  curve at $(0,0^-)$
 $$
 \tex{
 ``\,\,x^2= 4(1+\e)|t| \ln |\ln|t||\,\,", \quad \e>0,
 }
 $$
   he introduced the following explicit
  sub-solution
  {\em ``Barriere der
Irregularit\"at"}
  \cite[p.~58]{Pet34Cr}, \cite[p.~394]{Pet35} (see also comments in \cite{Moss00}):
  \be
  \label{sub1}
   \tex{
   ``\,\,\,v(x,t)= \frac {-1}{|\ln |t||^{1+\e_1}}\,\, {
   e}^{\frac {- x^2}{4t} \, k} + \frac 1{\ln|\ln |t||}\,\,\,"
   \quad \big(v_{\,t} < v_{\,xx}\big)
   \whereA
   \frac {1+\e_1}k < 1+\e,
    }
 \ee
and $\e_1<1$, $k<1$ are some positive constants. {\em Vice versa},
to prove regularity for
 $$
 \tex{
 ``\,\,\,x^2= 4|t| \ln |\ln|t||\,\,\,", \quad (\mbox{i.e., $\e=0$}),
 }
 $$
  the following explicit
  super-solution
  {\em ``Barriere der
Regularit\"at"}
  \cite[p.~58]{Pet34Cr} was used:
  \be
  \label{sub1NN} 
   \tex{
   ``\,\,\,v(x,t)= \frac {-1}{2|\ln |t||^{1\frac 16}}\,\, {
   e}^{\frac {- x^2}{4t}} + \frac 1{\sqrt[6]{|\ln |t||}} \,\,\,"
   \quad \big(v_{\,t} > v_{\,xx}\big).
    }
 \ee
In \cite{Pet35}, the standard nowadays notation $\underline
u(x,t)$ and $\overline u(x,t)$ for sub- and super-solutions were
used.
For more general boundaries given by
 \be
 \label{ssPP0}
 ``\,\,\,x^2 < 4t \ln \rho(t)\,\,\,",
 \ee
 Petrovskii \cite[p.~393]{Pet35} introduced the related function
 $\varphi(t)$ determined from
  \be
  \label{ssPP}
   \tex{
    ``\,\,\,\int\limits_{t_0}^t \frac{\rho(\eta) d \eta}\eta= 6 \log
    \varphi(t)\,\,\,",
 }
 \ee
 which was used in constructing suitable super- and sub-solutions.
Here, \ef{ssPP} is the origin of Petrovskii's criteria \ef{RR1}
given below.

In his another  paper \cite{Pet34} in 1934,
  Petrovskii also
  applied his novel method of {\em sub-} and {\em super-solutions}
  (i.e.,   before the classic Nagumo, or
  M\"uller--Nagumo--Westphal, Lemmas from the 1940-50s; for the Dirichlet problem, barrier-type
  ideas in the {\em m\'ethode de balayage}, or sweeping out, were used by Poincar\'e
  as early as in 1887 \cite{Poin87},
   with further development by Lebesgue in 1912 \cite{Leb12Bar} and by Perron in 1923
   for $\D u=0$ \cite[p.~432]{Pet34})
  in the study of
  randomness that had a great influence on random processes
  theory, which he connected with linear elliptic equations. This role was first explained in Khinchin's
   monographs in 1933 and 1936, \cite{Khin36}; see below.
Note that, in \cite{Pet34},
 Petrovskii dealt with the Dirichlet problems (``Dirichletschen
 Problems") for elliptic equations such as \cite[p.~428]{Pet34}
  $$
   \tex{
  ``\,\,\, M_{xx} \frac{\partial^2 u}{\partial x^2}+
  2 M_{xy} \frac{\partial^2 u}{\partial x \partial y}+
   M_{yy} \frac{\partial^2 u}{\partial y^2}+2  M_{x} \frac{\partial u}{\partial x}
   +2  M_{y} \frac{\partial u}{\partial y}=0\,\,\,",
   }
   $$
   where ``$M_{xx}M_{yy} \ge M_{xy}^2$", and, surely, Petrovskii's
   typical
  double $\ln \ln$ appeared as an its solution ``$u= \ln \ln \frac
  1{\sqrt{x^2+y^2}}$" on pp.~437--438.
   In three years time, sub-super-solution  techniques were applied
   in the seminal KPP-paper, 1937 \cite{KPP}.

\ssk

Petrovskii's 1934's paper
 \cite{Pet34Cr} is fully devoted to the optimal regularity and existence-uniqueness analysis
 for the heat equation, where he already introduced   his
 $2\sqrt{\ln\ln}$-criterion including the converging/diverging integral as in \ef{RR1} for
 irregularity/regularity. For instance, he showed that the
 ``5-$\log$ curve" and so on \cite[p.~57]{Pet34Cr} (cf.
 \cite[p.~404]{Pet35})
 $$
  \tex{
 ``\,\,\,x^2=4|t|[\ln|\ln t||+ 1 \frac 12 \ln \ln |\ln|t||+
 \ln \ln |\ln|t||
 (1+\e) \ln \ln
 \ln \ln |\ln |t||] \quad \mbox{u. s. w.}\,\,"
 }
 $$
is regular for $\e \le 0$ and is irregular for any $\e>0$. In
addition, for the 2D heat equation
 $$
 u_t= u_{xx}+u_{yy},
 $$
 it is mentioned \cite[p.~59]{Pet34Cr}  that
 the regularity/irregularity issues  occur about the
 surface
  $$
  ``\,\,\,x^2+y^2=-4(1+\e) t \ln|\ln|t||\,\,\,".
 $$
 In particular, then \cite[p.~397]{Pet35} the super-solution took
 the same form \ef{sub1}, where $x^2$ is replaced by $x^2+y^2$.
 In the later paper \cite{Pet35} in 1935, extensions to the heat equation in $\re^k \times \re_+$
were mentioned.

 Hence,
we  attribute Petrovskii's criterion to 1934 \cite{Pet34Cr}, and
not to 1935, \cite{Pet35}, as used had been done in most of other
``parabolic" regularity papers, we have referred to here.
  Concerning related earlier papers on parabolic PDEs  in
  non-cylindrical
  domains,  Solonnikov (1965)
  \cite{Sol65} and Ivasishen (1969) \cite{Iv69} should be mentioned.

\ssk

Thus, using such novel barriers,
 Petrovskii  \cite{Pet34Cr, Pet35} established the following:
 \be
 \label{PP1}
 \fbox{$
 \begin{matrix}
{\rm (i)}\,R(t)=2 \sqrt{-t} \,\,\sqrt{\ln|\ln (-t)|} \LongA
\mbox{$(0,0)$ is regular}, \andA \qquad\quad \ssk\\ {\rm (ii)}\,
R(t)=2(1+\e)\sqrt{-t}\,\, \sqrt{ \ln|\ln (-t)|}, \,\, \e>0 \LongA
 \mbox{$(0,0)$ is irregular}.
  \end{matrix}
  $}
  \ee
More precisely, he showed that, for the curve expressed in terms
of
 a
positive function $\rho(h) \to 0^+$ as $h \to 0^+$ ($\rho(h) \sim
\frac 1{|\ln h | }$ is about right) as follows:
 \be
 \label{RR1NN}
  \tex{
  R(t) = 2 \sqrt{-t} \,\,\sqrt{-\ln \rho(-t)},
   }
  \ee
 the {\em sharp regularity criterion} holds:
  \be
  \label{RR1}
   \fbox{$
\int\limits_0 \frac{{\rho(h)} \sqrt{|\ln \rho(h)|}}h \, {\mathrm
  d}h < \,(=)\, +\iy
   \LongA \mbox{$(0,0^-)$ is irregular (regular)}.
 $}
  \ee
 It is worth underlying again that both converging (irregularity)
 and diverging (regularity) integrals in \ef{RR1} as Dini--Osgood-type
 regularity criteria already appeared in the first Petrovskii paper
 \cite[p.~56]{Pet34Cr} of 1934. It should be mentioned that, since
 Petrovskii used  proposed him approach of constructing super- and
 sub-solutions, solving such partial differential inequalities
 led (as usual) unavoidably to a technical assumptions: for
 the super-solution (regularity of $(0,0)$), this is
 \cite[p.~392]{Pet35}
 $$
 ``\,\,\, 2. \quad t \log \rho(t) \to_{t \to 0}0; \,\,\,"
 $$
 and, for
 the sub-solution (irregularity of $(0,0)$)
 \cite[p.~397]{Pet35},
 $$
  \begin{matrix}
 ``\,\,\, 3. \quad t \log \rho(t) \to 0 \quad \mbox{f\"ur} \quad t \to 0;
 \,\,\,"\quad \mbox{and} \ssk\\
``\,\,\, 5. \quad \mbox{es ist} \quad \frac{\rho'(t)}{\rho(t)} \le
\frac
 1{8|t|}.\,\,\," \qquad\qquad\qquad\,\,\,
 \end{matrix}
 $$
 These purely technical assumptions can be got rid of, \cite{Abd00};
 our blow-up spectral-boundary layer approach is also
 assumptionless and rigorous for the heat equation, see Section \ref{S3.8}.

\ssk

The results of \cite{Pet34Cr, Pet35} were key important for
 probability theory, where this result is expressed as follows: if
 $x_t$ is a 1D Wiener process and $\Phi(t)>0$ is monotone such
 that
  \be
  \label{999}
 \tex{
  \int
  \limits
  ^\iy \frac{\Phi(t)}t\,\,{\mathrm e}^{-\Phi^2(t)/2}\, {\mathrm d}t=
  \iy,
  }
   \ee
  then, with unit probability, $|x_t| \le \Phi(t)$ for all $t \gg
  1$. On the contrary (this is about the criterion), if the
  integral converges, then with unit probability
  $\exists\,\,\{t_k\} \to \iy$ such that $|x_{t_k}| > \Phi(t_k)$;
  an
  optimal extension was  done in \cite{Khin36}.
  Eventually, these Petrovskii's results led to the so-called the {\em Law
  of Iterated Logarithms} (LIL) \cite[p.~392]{Durrett}, which was discovered by Khinchin and
  Kolmogorov. Earlier, in 1924 Khinchin proved that, for a sequence
   $\{X_i\}_{i \ge 1}$  of independent random
   variables with values $\pm 1$ with the probability $\frac 12$,
   there holds:
 \be
 \label{LLL1}
  S_n \equiv X_1+X_2+...+X_n \le \sqrt{(2+\e) n \ln \ln n} \,\,\,\,\, \mbox{for all}
   \,\,\, n \gg 1 \,\,\,(\forall \,\, \e>0)
   \ee
   with the probability one.
 This  sharp estimate improved earlier Hausdorf's inequality (1913)
  $S_n = o(n^{\frac 12+\e})$,   Hardy--Littlewood's
 one (1914) $S_n \le C \sqrt{n \ln n}\,$, and Steinhaus'
 inequality (1922)
 $\limsup S_n/\sqrt{2n \ln n} \le 1$.
 As a final step, in 1929, Kolmogorov proved the LIL for any bounded
 independent random variables not assuming that the summands were
 identically distributed; see more details in
  \cite[Ch.~7]{Durrett} and \cite{GrigKel00}.
   Clearly, \ef{PP1} has probabilistic roots in \ef{LLL1}.

A Wiener-type criterion, already established  by Khinchin
\cite{Khin36} (an earlier result was in 1933) in a probability
representation, for the boundary regularity for problems such as
\ef{ir1} and in $\ren$, was derived by Landis in 1969
\cite{Land69}
in terms of converging/diverging series of potentials  of
shrinking sets involved; see also \cite{Evans82} for a criterion
via thermal capacity also obtained along the lines of derivation
of Wiener's one for Laplace's equation.
  These results were not stated in the blow-up
evolution and more ``practical" Petrovskii's style \ef{RR1},
though Wiener's-type capacity criteria serve
  for more general types of
boundaries than according to Petrovskii's approach.

Petrovskii's integral criterion of the Dini--Osgood type given in
\ef{RR1}
    is true in the $N$-dimensional
radial case with (see \cite{Abd00}--\cite{Abd08Adv} for the recent
updating)
 $$
\sqrt{|\ln p(h)|} \quad \mbox{replaced by} \quad {|\ln
p(h)|}^{\frac N 2}.
 $$
 Several boundary regularity/irregularity results are now known for
a number of quasilinear parabolic equations
including degenerate porous medium operators. Though, some
difficult questions remain open even for the second-order
parabolic equations with order-preserving semigroups.
 We refer to \cite{Abd00, Abd05, Abd08Adv, Herr04, Mam01, Moss00} as a guide to a full history
and the extensions of these important results.

As far as we know, \ef{PP1} is the first clear appearance of the
``magic" $\sqrt{\rm log\,log}$ in PDE theory, currently associated
with the ``blow-up behaviour"  of the domain $Q_0$ and
corresponding solutions. Concerning other classes of nonlinear
PDEs generating blow-up $\sqrt{\rm log\,log}$ in other settings,
see references in \cite{GalLog}.



Finally, as an introduction to difficult features of higher-order
parabolic regularity to be studied, it is worth clearly stating
that
 \be
 \label{Pplus}
  \fbox{$
 \mbox{for solutions of changing sign of the HE, Petrovskii's
 criterion (\ref{RR1}) fails}
  $}
  \ee
  (obviously, positive subsolutions to prove irregularity are not applicable, while regularity
  can be still proved by comparison).
  In fact, a general criterion for arbitrary bounded changing sign data
 $u_0(x)$ in \ef{ir1} cannot be derived in principle. The origin
 of this is the same for the HE \ef{ir1}, the bi-harmonic one
 \ef{bih1}, and others; see our analysis below.

\subsection{Example: a refined asymptotics for the porous medium
equation}

In connection with the blow-up log-log, it is worth mentioning the
asymptotic result (seems, unique of this type) for the Dirichlet
problem as in \ef{ir1} \cite[\S~2]{Herr04}, which is now
formulated for the
 radial porous medium equation in the pressure form:
  \be
   \label{ss1}
    \begin{matrix}
    v_t= \D v^m, \quad m>1, \quad \mbox{where the pressure $u=v^{m-1}$
    solves:} \qquad\qquad\ssk\ssk\ssk\\
    u_t=u \D u + \g |\n u|^2 \inB Q_0, \quad \mbox{with}
    \quad R(t)=\sqrt{-t}, \quad \g= \frac 1{m-1}>0.\qquad\qquad
     \end{matrix}
     \ee
 Then the log-log occurs in the asymptotic behaviour of nonnegative
 solutions as $t \to 0^-$:
  \be
  \label{asss1}
   \tex{
  u(x,t) \sim \frac 1{4\ln|\ln(-t)|}\to 0
  \quad \mbox{uniformly on subsets}
  \quad \big| \frac {|x|}{\sqrt{-t}}-1 \big| \gg  \frac
  1{\ln|\ln(-t)|}.
   }
   \ee
 The computations leading to \ef{asss1} by asymptotic-matching
 ideas (some will be involved in our further study) are not easy,
  and a complete proof was not supplied in
 \cite{Herr04} (regardless the comment at the end of page 4 in
 \cite{Herr04}, it seems that a full and entirely rigorous
 justification of this kind of matching with ``floating" matching
 point can be extremely difficult and even seems illusive; this
 emphasizes a general complexity of such type of results even for
 second-order parabolic PDEs with the Maximum Principle).

 \ssk

 \noi{\bf Remark: on boundary criterion for {\bf \ef{ss1}}.}
 Obviously,  the delicate asymptotics \ef{asss1} allows one to
 detect the actual ``criterion" of the boundary regularity for the
 equation \ef{ss1} (strangely, this  was not addressed in
 \cite{Herr04}). We use the scaling invariance of \ef{ss1}:
  \be
  \label{kk1}
  u= A \, \hat u, \quad x= \sqrt A {\hat x }  \whereA A>0 \quad
  \mbox{is arbitrary}.
   \ee
   According to natural concepts of approximate similarity
   solutions (see e.g., \cite[Ch.~6]{SGKM}), we next assume that
   \be
   \label{kk2}
   A=A(t) \,\,\, \mbox{is not a constant, but a slow decaying
   function (relative to $\sqrt{-t}\,$).}
    \ee
  Then the equation for $\hat u(\hat x,t)$ will include an extra very
  small asymptotic perturbation, which will not affect the
  asymptotics. Finally, according to \ef{asss1} we set
   \be
   \label{kk3}
    \tex{
   A(t)=\frac 1{4\ln|\ln(-t)|}, \quad \mbox{so that}
   \quad \hat u \sim 1
   }
   \ee
   on the corresponding compact subsets as $t \to 0^-$.
 Thus, $u(0,0) \not = 0$ for the function
  \be
   \label{kk4}
   \hat R(t) \sim 2\sqrt {-t} \, \sqrt{\ln|\ln(-t)|},
  \ee
  and, in this irregularity condition, Petrovskii's  magic
  $\sqrt{\ln\, \ln}$ again mysteriously occurs.
  The constant ``2" here seems not relevant  as in \ef{PP1},
  since \ef{ss1} is nonlinear and admits also the symmetry
   $
   u \mapsto A
  u, \quad  t \mapsto \frac tA \quad \mbox{for any} \quad A>0,
  $
  so that
  the domain behaviour with
   $
  \hat R(t) \sim \frac 2{\sqrt{A}}\, \sqrt{-t}\,\sqrt{\ln|\ln(-t)|}
 $ is irregular for $(0,0)$.



\begin{center}

*\quad*\quad*

\end{center}


Thus, since the 1930s, Petrovskii's regularity
$\sqrt{\log\log}$-factor entered parabolic theory and generated
new types of asymptotic blow-up problems, which have been solved
for a wide class of parabolic equations with variable coefficients
as well as for some quasilinear ones. Nevertheless, such
asymptotic problems are very delicate and some of them  of
Petrovskii's type remain open even in the second-order case.

For higher-order poly-harmonic operators, the situation becomes
much more difficult, since
  the
Maximum and Comparison Principles and order-preserving properties
fail, so classic barrier techniques associated with sub-, super-,
and other solutions  are not applicable in principle.
 This
precisely reflects  the area of the present research: no Maximum
Principle tools, so that the boundary regularity problem in
Petrovskii's setting falls into the scope of a blow-up asymptotic
behaviour study.

 \section{Towards Petrovskii's criterion for the bi-harmonic
 flow: asymptotic blow-up setting, method, and layout of the paper}

\subsection{The basic Dirichlet problem under consideration}

 As a basic model, using the ``minimal" extension of the heat
 equation \ef{ir1} by increasing the order of the parabolic operator by two,
    we consider the {\em bi-harmonic equation} in the
 same shrinking domain as in \ef{ir1} with the zero Dirichlet
 boundary conditions on the lateral boundary $\partial Q_0$:
 \be
 \label{ir12}
  \left\{
 \begin{matrix}
   u_t=-u_{xxxx} \quad \mbox{in}
 \quad Q_0=\{|x|<R(t), \,\,\, -1 < t<0\},
 \ssk\ssk\\
  u=u_x=0 \atA x= \pm R(t), \quad -1 \le t <0,\quad\,\,\,
  \,\,\,\,\,\,\,
    \ssk\ssk\\
  u(x,0)=u_0(x) \onA [-R(-1),R(-1)],\qquad\quad\,\,\,\,\,\,\,\,
   \end{matrix}
    \right.
  \ee
where $u_0(x)$ is a  bounded and smooth function, $u_0(\pm
R(-1))=0$.

\ssk

First of all, obviously, in modern PDE theory, there are many very
strong regularity results, proved
for  wide classes of linear and nonlinear equations including, of
course, the simplest among others bi-harmonic one \ef{ir12}. We
have mentioned some key results and papers before and  will  not
indent to compete with those nice results and present no further
references to this classic literature.
 The only issue we propose here is as follows: to show how to derive
  \be
  \label{bl1}
 \mbox{a sharp boundary regularity criterion for (\ref{ir12})
 as a blow-up
 problem.}
  \ee
  As we have seen, these blow-up aspects of the regularity problem
 were already clearly revealed  by Petrovskii in 1934. Though,
 for the second-order case, there are other equivalent  approaches based on the positivity
 of the kernel and various order-preserving features, which fail for oscillatory kernels.

\ssk

To be more precise, this blow-up study for \ef{ir12} is performed
for:
 \be
  \label{gen}
  \mbox{a {\em generic} class of solutions exhibiting
 a ``centre subspace" behaviour.}
  \ee
 In other words, as our Hermitian spectral theory of higher-order
 non self-adjoint linear operators in Section \ref{S2} shows,
 one can expect that, in addition to \ef{gen},
 \be
  \label{genNN}
  \begin{matrix}
 \exists\,\,\,\, \mbox{an infinite-dimensional set of  ``stable
 manifold" patterns,}\ssk\\
 \mbox{for which the conditions of regularity of $(0,0)$ are all
 different.}
  \end{matrix}
 \ee

It is principal for us that,
 to confirm the optimal character of the class of generic blow-up
 asymptotics, we will check (Section \ref{S3.8}) that, by elementary spectral calculus,
 \be
  \label{genPet}
   \begin{matrix}
  \mbox{the same asymptotic approach to
 the classic problem (\ref{ir1}) by classic}\qquad\quad
 \ssk\\
  \mbox{Hermite spectral theory yields precisely Petrovskii's criterion in
 (\ref{RR1}).}\qquad\quad
 \end{matrix}
  \ee
Note that the problem: \ef{gen} or \ef{genNN}, also exists for the
heat equation \ef{ir1} (cf. \ef{Pplus}), but here, by the Maximum
Principle, positive solutions always belong to the generic class
\ef{gen}. For the bi-harmonic PDE \ef{ir12}, in a general setting,
distinguishing generic patterns \ef{gen} from those in \ef{genNN}
is difficult and seems even impossible (or makes no sense).


 \subsection{Slow growing factor $\var(\t)$}

 Thus,
similar to \ef{PP1}, we need to assume that
 \be
 \label{ph1}
 R(t)= (-t)^{\frac 14}\, \var(\t) \whereA \t= - \ln(-t) \to + \iy
 \asA t \to 0^-.
  \ee
 Here, the scaling main factor  $(-t)^{ 1/4}$ naturally comes from the bi-harmonic
 kernel variables (see \ef{ph3} and \ef{Fund}),
and $\var(\t)>0$ is an unknown  slow growing function
 satisfying\footnote{As we have mentioned, for the simpler case
 $\var(\t),\, \var'(\t) \to 0$ as $\t \to \iy$, the regularity was already proved
 by Miha${\rm\check{i}}$lov in 1963
 \cite{Mih63I, Mih63II}; in a certain sense, this extended the Gevrey-like result
 \ef{Gev1} for $m=1$; see \ef{Gev1Mih}.}
 \be
 \label{vv1}
  \tex{
 \var(\t) \to +\iy, \quad \var'(\t) \to 0, \andA \frac {\var'(\t)}{\var(\t)} \to 0 \asA \t \to
 +\iy.
  }
  \ee
 Moreover, as a sharper characterization of the above class of
 {\em slow growing functions}, we use the following criterion:
  \be
  \label{vv2}
   \tex{
   \big( \frac{\var(\t)}{\var'(\t)}\big)' \to \iy \asA \t \to +\iy
    \quad (\var'(\t) \not = 0).
    }
    \ee
    This is a typical condition in blow-up analysis distinguishing classes of
    exponential (the limit is 0), power-like (a constant $\not =
    0$), and slow-growing functions. See \cite[pp.~390-400]{SGKM},
    where in Lemma 1 on p.~400, extra properties of slow-growing
    functions \ef{vv2} are proved. For instance, one can derive
    the following comparison of such $\var(\t)$ with any power:
     \be
     \label{al1}
     \mbox{for any $\a>0$}, \quad \var(\t) \ll \t^\a \andA
     \var'(\t) \ll \t^{\a-1} \forA \t \gg 1.
     \ee
 Such estimates are useful
 in  evaluating
  perturbation terms in the rescaled equations.

Thus, the monotone positive function $\var(\t)$ in \ef{ph1}
  is assumed to determine a sharp behaviour of the
boundary of $Q_0$ near the shrinking point $(0,0)$ to guarantee
its regularity. In Petrovskii's criterion \ef{PP1}, the almost
optimal
 function, satisfying  \ef{vv1}, \ef{vv2}, is
  \be
  \label{ph2}
   \var_1(\t) \sim 2\sqrt {\ln \t} \asA \t \to + \iy,
   \ee
a dependence we have to compare  our final results for the
bi-harmonic equation with.



\subsection{First kernel scaling and layout}

 By \ef{ph1}, we
perform the similarity scaling
 \be
 \label{ph3}
  \tex{
  u(x,t)=  v(y,\t) \whereA y=\frac x{(-t)^{
  1/4}}.
  }
   \ee
   Then the rescaled function $v(y,\t)$ now solves the
   rescaled equation
 \be
 \label{ph4}
  \left\{
  \begin{matrix}
  v_\t= \BB^* v \equiv - v_{yyyy}- \frac 14 \, y v_y
  \inB Q_0=\{|y| < \var(\t), \,\,\, \t>0\}, \ssk\ssk\\
  v=v_y=0 \atA y= \pm \var(\t), \,\,\t \ge 0,
  \qquad\qquad\qquad\qquad\qquad\,
  \ssk\ssk\\
  v(0,y)=v_0(y) \equiv u_0(y) \onA
  [-R(-1),R(-1)].\qquad\quad\,\,\,\,\,\,\,\,\,\,\,\,
   \end{matrix}
   \right.
  \ee

 In view of the divergence \ef{vv1}, it follows that our final analysis will
 essentially depend on the spectral properties of the linear
 operator $\BB^*$ on the whole  line $\re$. We reflect this Hermitian spectral theory
in the next Section \ref{S2}, with application to the regularity
criterion in Section \ref{S3}. In particular, this  differs our
analysis from Kondrat'ev's classic one \cite{Kond66}.

However, in Section \ref{S.l}, we begin the study of the
regularity of the vertex $(0,0)$ of the bi-harmonic  {\em backward
fundamental parabolae}:
 \be
 \label{P1}
 R(t)=l(-t)^{\frac 14}, \quad \mbox{i.e.,} \quad \var(\t) \equiv
 l={\rm const.}>0.
  \ee
Then the problem \ef{ph4} is considered on the fixed interval
$I_l=\{|y|<l\}$, so that the final conclusion entirely depends on
spectral properties of $\BB^*$ in $I_l$ with Dirichlet  boundary
conditions. Of course, then the spectral problem for $\BB^*$ (not
a pencil) becomes a very particular case of the general setting
developed in \cite{Kond66}, but nevertheless, the clear conclusion
on regilarity/irregularity becomes rather involved, where numerics
are necessary to fix final details. In addition,
 as we pointed out,  in more general setting for the
fundamental backward paraboloids  in $\ren$, the existence,
uniqueness, and regularity of solutions in Sobolev spaces  was
proved in a number of papers such as \cite{Mih61, Mih63I, Mih63II,
Fei71}, etc.
 Note that in \cite[p.~45]{Mih63I}, the zero boundary data was
 understood in the {\em mean sense} (i.e., in the $L^2$-sense along a sequence smooth internal contours
 ``converging" to the boundary).
 Nevertheless,
we have to stress attention to this simple case in order to reveal
the exact transition between regularity and irregularity in the
classic sense in the critical case \ef{P1}. It seems that such a
border case in between was not addressed before sharply enough.

 Thus, our conclusion is rather disappointing: unlike the
classic heat equation \ef{ir1}, for the bi-harmonic equation
\ef{ir12},  {\em the vertexes   of the fundamental parabolae
\ef{P1} are not necessarily
 regular}\footnote{It is worth comparing this with the following
 well-known (since 1986)
 negative conclusion in elliptic theory: for fourth-order elliptic
 equations (with $N \ge 8$; for $5 \le N \le 7$, the biharmonic capacity settles the regularity
 result in Wiener's sense \cite{Maz79}; $N<4=2m$ is done by Sobolev embedding),
  the vertex of a cone can be irregular-singular (the solution unbounded at the vertex),
 \cite{MN86}.}. For instance, on the basis of careful numerics, we show that:
 \be
 \label{P2}
 \fbox{$
 \begin{matrix}
 R(t)=4(-t)^{\frac 14} \LongA
 \mbox{$(0,0)$ is regular, while}
  \ssk\ssk\\
 R(t)=5(-t)^{\frac 14} \LongA \mbox{$(0,0)$ is irregular.}\qquad
 \end{matrix}
 $}
 \ee
 Moreover, in the latter case, the vertex $(0,0)$ is {\em
 singular}, i.e.,
 \be
 \label{P99}
 \mbox{in general, $u(x,t)$ is unbounded as $(x,t) \to (0,0^-)$}.
  \ee
We also claim that (by a continuity argument) there exists some
$l_1 \in (4,5)$ such that $(0,0)$ is irregular, but not singular,
i.e., the limit $u(0,0^-)$ exists.

These regularity variations for constant functions $\varphi(\t)
\equiv l$ inevitably impose a special restriction for the study in
Section \ref{S3} unbounded $\var(\t)$'s satisfying \ef{vv1}, where
we derive a Petrovskii-like ``criterion" of regularity. Namely,
the irregularity condition is obtained in terms of an
Osgood--Dini-type integral condition, while, in an accordance with
\ef{P2}, for the regularity, a special ``oscillatory cut-off" of
the lateral boundary must be applied.

\ssk

In Section \ref{S4}, we discuss the extensions of the boundary
regularity analysis to  $2m$th-order poly-harmonic operators
 \be
 \label{pol1}
 u_t= (-1)^{m+1} D_x^{2m}u \inB Q_0 \quad (m \ge 2),
  \ee
  with the zero Dirichlet conditions, as well as for the same
$N$-dimensional  problem  for \ef{PE1} in $Q_0$,
where the shape of the shrinking domain $Q_0$ is also of
importance for the boundary regularity. For instance, one can
consider the radial $Q_0$ given in \ef{ir12}, where $|x|$
 stands for the radial spatial
 variable. We also discuss regularity conditions for the
third-order linear dispersion equations, for a quasilinear
fourth-order porous-medium-type equation (the PME--4), and for the
linear wave (beam) equation of the fourth order.

\ssk

The present paper aims to give a first insight into principal
difficulties of a  sharp study of boundary regularity for
higher-order parabolic and other evolution PDEs. Such a study
inevitably generates a number of difficult mathematical problems,
which, partially, do not appear in the second-order case, or can
be avoided by positivity kernel properties.
 We must admit that
some of them
 are not and even cannot be completely rigorously
solved. The main problem of concern is that distinguishing generic
patterns in  \ef{gen} from non-generic ones in \ef{genNN} is not
possible in general. Here, delicate boundary layer theory is
essentially involved, where
  we ought to use some accurate numerical calculations by the
{\tt MatLab} whenever necessary to avoid huge technical
digressions (using those seems also inevitable for truly
$2m$th-order parabolic equations with large enough $m \ge 2$).
 Regardless a  certain rigorous  incompleteness of the mathematical analysis in
 some steps (which seems to be inevitable in general), we decide to demonstrate the whole machinery of
 the asymptotic blow-up methods in parabolic boundary regularity theory, and
 hope that this will help to attract some attentive Readers to
  improve the results when necessary.
 It is worth mentioning that, as can be seen in many directions of
 modern PDE theory of linear and nonlinear equations, the
 transition to higher-order models is accompanying by a dramatic
 increase of the complexity of the methods involved to achieve in
 a rigorous manner the necessary (basic or not) results. This
 increase of complexity sometimes measures in orders, and often
  desired rigorous proofs can be  illusive in a sufficient generality.
 A suitable restriction of the generality to achieve  rigorous
 conclusions, though being an important research step, can lead to
 an extended and often unjustified number of artificial hypothesis
 and even yield  non-constructive assumptions (i.e., those that cannot be
 checked in a reasonable finite time or at all).
 In what follows, we prefer to avoid stating such theorems with
 non-constructive hypothesis, though this is indeed plausible in a
 few places\footnote{``The main goal of a mathematician is not
 proving a theorem, but an effective investigation of the
 problem...", A.N.~Kolmogorov, 1980s (the author apologizes for a
non-literal translation from the Russian).}.
 In this case,
we believe that  the exchange of the ideas and methods, even in
the case of
  a certain lack of a completely rigorous justification,
 can be key for further improvements and  developing more consistent mathematical PDE
 theory. Of course, this imposes certain restrictions on the paper
 style, but the author hopes that the interested Reader will
 easily distinguish the rigorous and non-rigorous arguments to be
 used.

\section{Fundamental solution and Hermitian spectral theory for
$\{\BB, \,\BB^*\}$}
 \label{S2}

For convenience of the further study, we consider the $2m$th-order
poly-harmonic equation \ef{pol1}, and will describe the necessary
 spectral properties of the linear differential operator (the
 analogy of that in \ef{ph4} for any $m \ge 2$)
\begin{equation}
 \label{B1*}
  \tex{
 \BB^* = (-1)^{m+1} D^{2m}_y - \frac {1}{2m}
 \, y D_y,
 }
 \end{equation}
and of its adjoint $\BB$ in the standard $L^2$-metric given by
  \begin{equation}
 \label{B1}
  \tex{
 \BB = (-1)^{m+1} D_y^{2m} + \frac {1}{2m}
 \, y D_y +  \frac{1}{2m}\, I \quad (I \,\,\, \mbox{denotes the identity}).
 }
 \end{equation}
 Both operators are not symmetric and do not admit a self-adjoint
extension. We will follow \cite{Eg4} in presenting necessary
spectral theory.

\subsection{The fundamental solution and its sharp estimates}

We begin with  determining the spectrum and the eigenfunctions of
the adjoint operator $\BB$, which appears in constructing  the
{\em fundamental solution} of \ef{pol1} that takes
 the standard
self-similar form
 \begin{equation}
\label{Fund}
 \tex{
 b(x,t) = t^{-\frac 1{2m}}F(y), \quad y= \frac x{t^{1/2m}}.
 }
\end{equation}
Substituting $b(x,t)$ into (\ref{pol1}) yields that the radially
symmetric profile
 $F(y)$ is the unique even  integrable solution of the linear ODE
 \begin{equation}
\label{ODEf}
 \tex{
 \BB F  = 0
 \quad {\rm in} \quad \re, \quad \int_\re F =1,
 }
\end{equation}
so it is a null eigenfunction of $\BB$. Taking the Fourier
transform leads to
 \be
 \label{FundSol}
  \tex{
 F(y) = \a_0  \int\limits_0^\infty {\mathrm e}^{-s^{2m}} \cos(s y)\, {\mathrm
 d}s,
 }
 \ee
 where $\a_0>0$ is the normalization constant,
 and, more precisely \cite{Eid},
\begin{equation}
\label{Eidf}
 \tex{
  F(y) = \frac 1 {\sqrt{2\pi}}\,  \int\limits_0^\infty {\mathrm
  e}^{-s^{2m}}
\sqrt{s|y|}\,\, J_{- \frac 1 2}(s|y|)\, {\mathrm d}s \quad
\mbox{in} \,\,\, \re, }
\end{equation}
where $J_\nu$ denotes Bessel's function.
   The rescaled kernel $F(y)$ then
satisfies a standard pointwise exponential estimate \cite{EidSys}
 \be
 \label{es11}
  \tex{
  |F(y)| \le  D \,{\mathrm e}^{-d_0|y|^{\a}}\quad \mbox{in} \,\,\, \re
  \whereA \a = \frac{2m}{2m-1}
  }
   \ee
 and $D$ and $d_0$ are some positive constants ($d_0$ to be specified below).
 Such optimal exponential estimates of the fundamental solutions
 of higher-order parabolic equations are well-known and were first
 obtained by Evgrafov--Postnikov (1970) and Tintarev (1982); see
 Barbatis \cite{Barb, Barb04} for key references and results.

However, as a crucial  issue for the further boundary point
regularity study, we will need a sharper, than given by \ef{es11},
asymptotic behaviour of the rescaled kernel $F(y)$ as $y \to
+\iy$. To get that, we re-write the equation \ef{ODEf} on
integration once as
 \be
 \label{i1}
 \tex{
 (-1)^{m+1}F^{(2m-1)} + \frac 1{2m} \, y F=0 \inB \re.
 }
 \ee
 Using standard classic WKBJ-type asymptotics, we substitute into \ef{i1}
 the function
  \be
  \label{i2}
  F(y) \sim y^{-\d_0} \, {\mathrm e}^{a y^\a} \asA y \to + \iy,
   \ee
   exhibiting two scales.
This gives the algebraic equation for $a$,
 \be
 \label{i3}
  \tex{
 (-1)^m (\a a)^{2m-1}= \frac 1{2m}, \andA \d_0=  \frac{m-1}{2m-1}>0\,
 .
 }
  \ee
 Note that the slow algebraically decaying factor $y^{-\d_0}$ in
\ef{i2} is available for any $m \ge 2$ and is absent for $m=1$ for
the pure exponential positive Gaussian profile $F$; see below.

By construction, one needs to get the root $a$ of \ef{i3} with the
maximal ${\rm Re}\, a<0$. This yields (see e.g., \cite{Barb,
Barb04} and \cite[p.~141]{GSVR})
 \be
 \label{i4}
  \tex{
 a= \frac{2m-1}{(2m)^\a} \big[ \cos\big( \frac{m \pi}{2m-1}\big) +
 \ii \sin\big( \frac{m \pi}{2m-1}\big)\big] \equiv -d_0 + \ii b_0
 \quad \big(\,d_0 \equiv \sin \big( \frac \pi{2(2m-1)}\big)>0\, \big).
 }
 \ee
Finally, this gives the following double-scale asymptotic of the
kernel:
 \be
 \label{i5}
  \tex{
  F(y) =
   y^{-\d_0} \, {\mathrm e}^{-d_0 y^\a} \big[ C_1 \sin (b_0 y^\a)+
   C_2 \cos (b_0 y^\a)\big]+... \asA y \to + \iy ,
   }
   \ee
 where $C_{1,2}$ are real constants, $|C_1|+|C_2| \not = 0$.
 In \ef{i5}, we present the first two leading terms from the
 $m$-dimensional bundle of exponentially decaying asymptotics.

In particular, for the bi-harmonic equation \ef{ir12}, we have
 \be
 \label{i6}
  \tex{
  m=2: \quad \a= \frac 43,  \quad d_0=3 \cdot 2^{-\frac{11}3},
 \quad b_0=3^{\frac 32} \cdot 2^{-\frac{11}3},
  \andA \d_0= \frac
  13.
   }
   \ee

\subsection{The discrete real spectrum and eigenfunctions of
 $\BB$}

We describe the spectrum  $\sigma(\BB)$ of $\BB$
 in
the space $L^2_{\rho}({\bf R})$  with the exponential weight
  \begin{equation}
 \label{rho11}
  \tex{
 \rho(y) = {\mathrm e}^{a |y|^{\a}}>0 \quad {\rm in} \,\,\, \re \quad
 \big(\a=\frac{2m}{2m-1}\big),
 }
 \end{equation}
where $a \le 2d_0$ is a  positive constant. Denoting by $\langle
\cdot, \cdot \rangle_\rho$ and $\|\cdot\|_\rho$ the corresponding
inner product  and the induced norm respectively, we introduce a
standard  Hilbert (a weighted Sobolev) space of functions
$H^{2m}_{\rho}(\re)$ with the inner product and the norm
\[
 \tex{
 \langle v,w \rangle_{\rho} = \int\limits_{\re} \rho(y) \sum\limits_{k=0}^{2m}
 D^{k}_y
 v(y) \, \overline {D^{k}_y w(y)} \,{\mathrm d} y, \quad
\|v\|^2_{\rho} = \int\limits_{\re} \rho(y) \sum\limits_{k=0}^{2m}
|D^{k}_y
 v(y)|^2 \, {\mathrm d} y.
  }
\]
Then $H^{2m}_{\rho}(\re) \subset L^2_{\rho}(\re) \subset
L^2(\re)$, and  $\BB$ is a bounded linear operator from $
H^{2m}_{\rho}(\re)$ to $ L^2_{\rho}(\re)$. With these definitions,
the spectral properties of the operator $\BB$ are given by:

\vspace{0.1in}

\begin{lemma}
\label{lemspec}
 {\rm (i)}  The spectrum of $\BB$
comprises real simple eigenvalues only,
 \begin{equation}
\label{spec1}
 \tex{
 \sigma(\BB)=
\big\{\lambda_k = -\frac k{2m}, \,\, k = 0,1,2,...\big\}.
 }
\end{equation}
 {\rm (ii)} The eigenfunctions $\psi_k(y)$ are given by
 \begin{equation}
\label{eigen}
 \tex{
 \psi_k(y) = \frac{(-1)^{k}}  {\sqrt{k !}}\,
  D^k_y F(y)
  }
\end{equation}
and form a complete set
 in $L^2({\re})$ and in $L^2_{\rho}({\re})$.

\noi {\rm (iii)} The resolvent $(\BB-\lambda I)^{-1}$
for  $\lambda \not \in \sigma(\BB)$ is a compact integral operator
in $L^2_{\rho}(\re)$.
\end{lemma}


 The operators $\BB$ ($\BB^*$) have zero Morse index (no
eigenvalues have positive real part).

\subsection{The polynomial eigenfunctions of the operator $\BB^*$}

We now consider the operator (\ref{B1*}) in the weighted space
$L^2_{\rho^*}(\re)$ ($\langle \cdot, \cdot \rangle_{\rho^*}$ and
$\|\cdot\|_{\rho^*}$ are the inner product and the norm)
 with the exponentially decaying weight
function
  \begin{equation}
\label{rho2}
 \tex{
 \rho^*(y) \equiv \frac  1 {\rho(y)} = {\mathrm e}^{-a|y|^{\a}} > 0,
 }
\end{equation}
 and ascribe to $\BB^*$ the domain
 $H^{2m}_{\rho^*}(\re)$, which is dense in $L^2_{\rho^*}(\re)$.
 Then
$
 \BB^*: \,\, H^{2m}_{\rho^*}(\re) \to L^2_{\rho^*}(\re)
$
 is a bounded linear operator. Hence, $\BB$ is adjoint
 to $\BB^*$ in the usual sense: denoting by $\langle \cdot,\cdot \rangle $  the
inner product on $L^2(\re)$, we have
  \begin{equation}
 \label{Badj1}
 \langle \BB v, w \rangle =  \langle v, \BB^* w \rangle
 \quad \mbox{for any} \,\,\, v \in H^{2m}_\rho(\re), \,\,\,
 w \in H^{2m}_{\rho^*}(\re).
 \end{equation}
The eigenfunctions of $\BB^*$ take a particularly simple
polynomial form and are as follows:


 \begin{lemma}
\label{lemSpec2}
 {\rm (i)} $ \sigma(\BB^*)=\s(\BB)$.

 \noi{\rm (ii)} The eigenfunctions  $\psi^*_k(y)$ of $\BB^*$ are polynomials
in $y$ of order $k$ given by
  \begin{equation}
 \label{psi**1}
  \tex{
 \psi^*_k(y) = \frac 1{\sqrt{k !}}
\sum\limits_{j=0}^{[k/2m]} \frac {(-1)^{m j}}{j !}D^{2m j}_y y^k,
\quad k=0,1,2,...\,,
 }
 \end{equation}
and form a complete subset  in $L^2_{\rho^*}(\re)$.

  \noi {\rm (iii)}
$\BB^*$ has compact resolvent $(\BB^*-\lambda I)^{-1}$ in
$L^2_{\rho^*}(\re)$ for $\lambda \not \in \sigma(\BB^*)$.
\end{lemma}


 With the definition (\ref{psi**1}) of the adjoint basis, by
integrating by parts,  the orthonormality condition holds
 ($\d_{kl}$ is the Kronecker delta):
\begin{equation}
 \label{Ort}
\langle \psi_k, \psi_l^* \rangle = \d_{kl} \quad \mbox{for any}
\,\,\,\, k,\,\,l \ge 0.
 \end{equation}
  For $m=2$ (this case will be treated in greater
  detail),
  the first eigenfunctions are
 \be
 \label{psi44}
  \begin{matrix}
 \psi_0(y) = 1, \quad \psi_1(y)=y, \quad \psi_2(y) = \frac 1{\sqrt{2}}\,  y^2,
 \quad \psi_3(y) = \frac 1{\sqrt 6}\, y^3,\qquad \ssk\ssk\\
  \psi_4(y)
=  \frac 1 {\sqrt{24}}\, (y^4 + 24),
 \,\, \psi_5(y)= \frac 1{2 \sqrt{30}}\,(y^5+ 120\, y), \,\,
 \psi_6(y) =\frac 1 {\sqrt{6!}}\, (y^6 +360y^2), \qquad
 \end{matrix}
 \ee
 etc.,
  with the corresponding eigenvalues $0$, $-\frac 14 $, $-\frac 12$, $-\frac 34$, $-1$, $- \frac 54$,
    $- \frac 32$.


%

\section{The vertex of fundamental parabolae can be regular or
irregular}
 \label{S.l}

Thus, consider the backward  fundamental parabolae given by
\ef{P1}. We continue denote by $\BB^*$ the corresponding linear
operator
 \be
 \label{P3}
  \tex{
 \BB^*v=-v^{(4)}- \frac 14\, y v' \inB I_l, \quad v=v'=0 \atA y=
 \pm l,
 }
 \ee
 with the standard definition of the domain, etc.
 Here,
$\BB^*$ is a {\em regular} ordinary differential operator with
bounded smooth coefficients and a discrete spectrum $\s(\BB^*)$
\cite{Nai1}. Note that $\BB^*$ is not symmetric and does not admit
a self-adjoint extension, so that $\s(\BB^*)$ is not necessarily
real (though we will observe a lot of real eigenvalues to be
explained by a branching phenomenon).
 By $\l_k(l)$, we denote $l$-branches of eigenvalues of $\BB^*$
 depending on the length $l>0$. In view of Section \ref{S2}, the limit $l \to + \iy$
is of a particular interest.

The regularity criterion of the vertex  is now easy:
 \be
 \label{P4}
  \fbox{$
  \mbox{vertex $(0,0)$ of  parabolae (\ref{P1}) is regular
  iff  $\s(\BB^*) \subset {\mathbb C}_- = {\mathbb
  C}\cap\{{\rm
  Re}\, \l <0\}$}.
  $}
  \ee
  Then the solution $v(y,\t)$ (and any $u(x,t)$) gets exponentially small: as $\t \to
  +\iy$,
 \be
 \label{P41}
 v(y,\t) \sim {\mathrm e}^{\l_0 \t}\to 0 \LongA
 |u(x,t)| \sim  O\big(t^{- {\rm Re}\, \l_0}\big)\to 0 \asA t \to 0^-,
 \ee
 where $\l_0(l)$ is the eigenvalue of $\BB^*$ with the maximal negative ${\rm Re}$ (in fact,
 $\l_0(l)$ turns out to be  always real).
Let us now briefly describe main steps of the current study.

\subsection{The vertex $(0,0)$ is regular for not that large
$l < 3.9779...\,$}

This conclusion is simple: as a standard practice, by multiplying
 in $L^2(I_l)$ the eigenfunction equation
 \be
 \label{P5}
  \tex{
 \BB^* \Psi^* \equiv - (\Psi^*)^{(4)}- \frac 14\, y (\Psi^*)' = \l
 \Psi^*
  }
  \ee
  by the complex conjugate $\bar \Psi^*$ and summing up with the
  result of multiplying the complex conjugate equation by
  $\Psi^*$, on integration by parts, one obtains
 \be
 \label{P6}
  \tex{
  (\l + \bar \l) \int|\Psi^*|^2=-2 \int |\Psi''|^2 + \frac 14
  \int|\Psi^*|^2.
  }
  \ee
 Using the Poincar\'e inequality:
  \be
  \label{P7}
   \tex{
   \int |\Psi''|^2 \ge \Lambda_0(l) \int |\Psi|^2, \quad \Psi \in
   H^2_0(I_l); \quad \Lambda_0(l)= \frac 1{l^4}\, \Lambda_0(1),
   \quad \Lambda_0(1)=31.3...\, ,
    }
    \ee
      where $\Lambda_0(1)$ is the first eigenvalue of $D^{(4)}_y>0$ in
$H^2_0(I_1)$, we obtain from \ef{P6}:
 \be
 \label{P8}
 \tex{
(\l + \bar \l) \int|\Psi^*|^2 \le - \big[ \frac{2
\Lambda_0(1)}{l^4} - \frac 14 \big] \int |\Psi^*|^2.
  }
   \ee
   Therefore, the regularity of the vertex $(0,0)$ is guaranteed for
   (this is not sharp but close):
 \be
 \label{P9}
 \tex{
 \frac{2
\Lambda_0(1)}{l^4} - \frac 14>0 \LongA l<l_*=[8
\Lambda_0(1)]^{\frac 14}= 3.9779...\,.
 }
 \ee

\subsection{First eigenvalues of $\BB^*$ are real}

Here and later on in a few special technically difficult cases, we
are going to rely on a clear numerical evidence for our
conclusions. To this end, we use the {\tt MatLab} solver {\tt
bvp4c} with the enhanced accuracy and tolerances
 \be
 \label{Tol1}
 \mbox{
from $10^{-10}$ to the minimal admitted $10^{-13}$.}
  \ee
 This will
guarantee that real eigenvalues, which can be very small, are
computed correctly.

\begin{figure}
\centering
\includegraphics[scale=0.65]{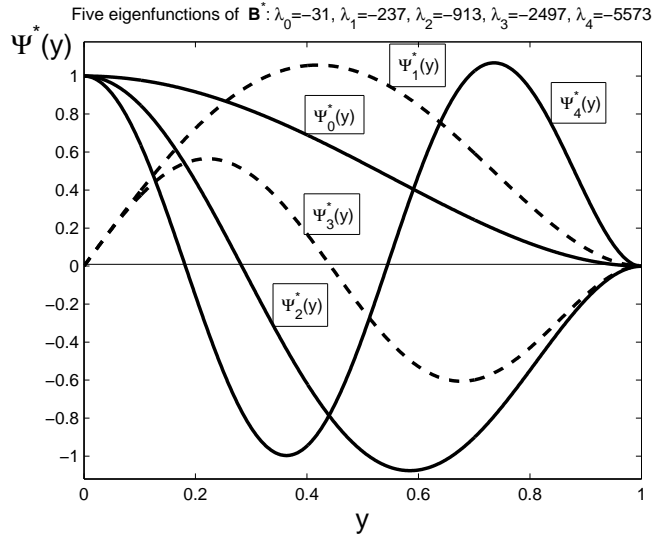} 
\vskip -.3cm \caption{\small Five first real eigenfunctions of the
operator \ef{P3} for $l=1$.}
   \vskip -.3cm
 \label{FP1}
\end{figure}

In Figure \ref{FP1}, we demonstrate first eigenfunctions of the
operator $\BB^*$ in \ef{P3} for $l=1$, where $\Psi_0^*$,
$\Psi_2^*$, and
 $\Psi_4^*$ are even functions, while $\Psi_1^*$ and $\Psi_3^*$ are odd ones.
Note that, by obvious reasons, the first eigenvalue
 $$
 \l_0(1)=-31.16...
 $$
is rather close to the eigenvalue $-\Lambda_0(1)=-31.3...$ in
\ef{P7} of the self-adjoint counterpart $-D^{(4)}_y$, i.e., the
non-symmetric perturbation $- \frac 14 \, y D_y$ is negligible.
Actually, the ``Sturmian geometric structure" of eigenfunctions in
Figure \ref{FP1}: each $\Psi_k^*(y)$ has precisely $k$ zeros in
$I_1$, which is true for $-D^{(4)}_y$ \cite{Ellias}, clearly
remains valid for $\BB^*$. Moreover, we have observed that {\em
the first eigenfunctions of $\BB^*$ and $-D^{(4)}_y$ for $l=1$
coincide within the accuracy $0.5 \times 10^{-2}$.} In fact, this
suggests to get $\l_k(l)$ by branching at $\mu=0^+$ from
eigenfunctions of $-D^{(4)}_y<0$ by using the following family for
$\mu \in \big[0, \frac 14\big]$
 (a homotopy path) of operators in $H_0^4(I_l)$:
 \be
 \label{771}
  \tex{
  \BB^*_\mu=- D^{(4)}_y -\mu \, y D_y' = \BB^* \forA \mu =
  \frac 14.
  }
  \ee
For $\mu=0$, the self-adjoint $\BB^*_0= -D^{(4)}_y$ has the
desired real spectrum and complete-closed set of eigenfunctions;
see \cite{Kato} and \cite{VainbergTr} for classic
perturbation/branching theory.
 Later on, we
present another branching explanation of the origin of real
eigenvalues/eigenfunctions of the non self-adjoint operator
$\BB^*$ by using results from Section \ref{S2} reflected $l=\iy$.

\subsection{$\l_0(l)$ changes sign: numerical evidence}

\begin{figure}
\centering \subfigure[first zero $l_1=4.08...$]{
\includegraphics[scale=0.52]{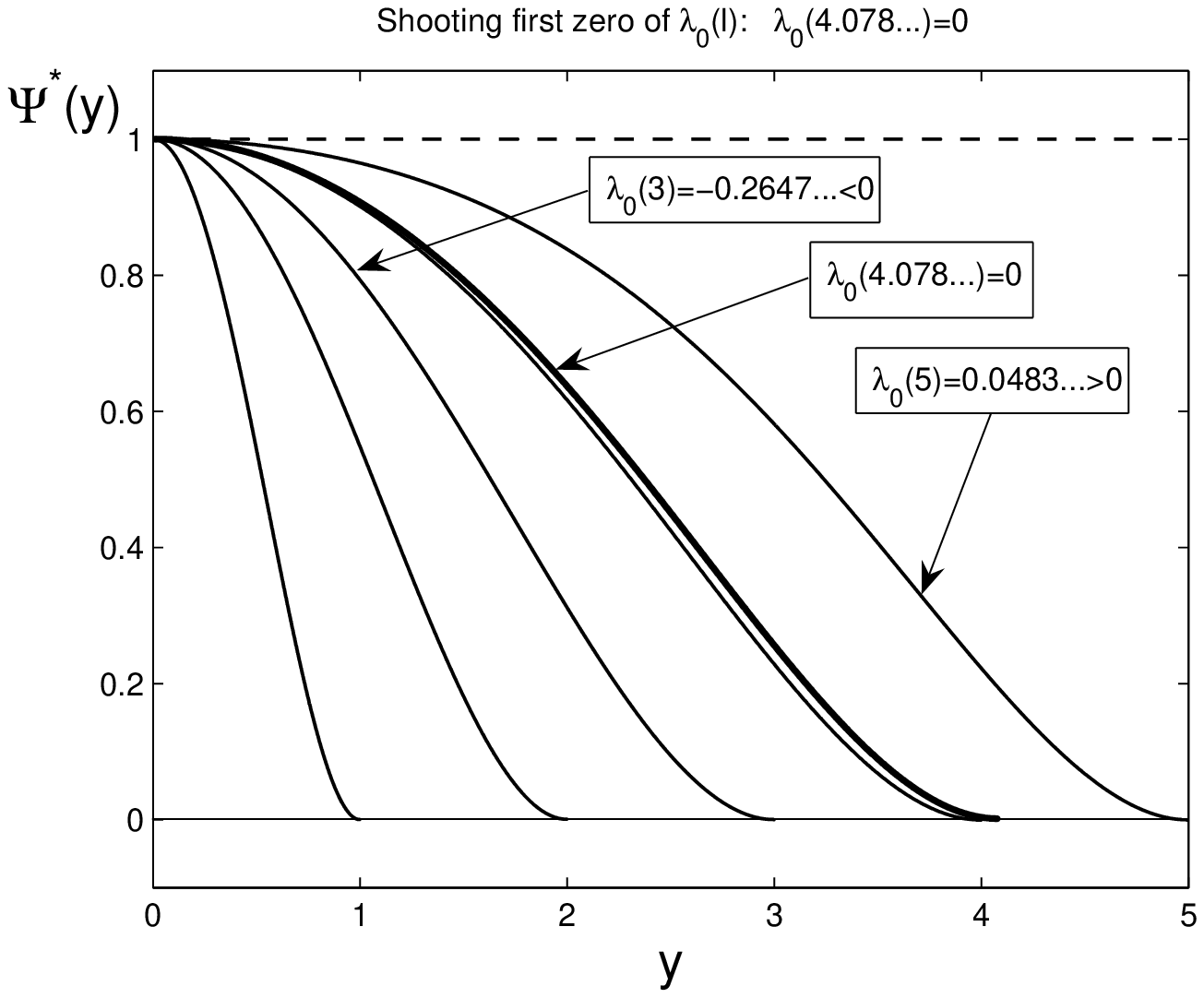}
} \subfigure[second zero $l_2=7.25...$]{
\includegraphics[scale=0.52]{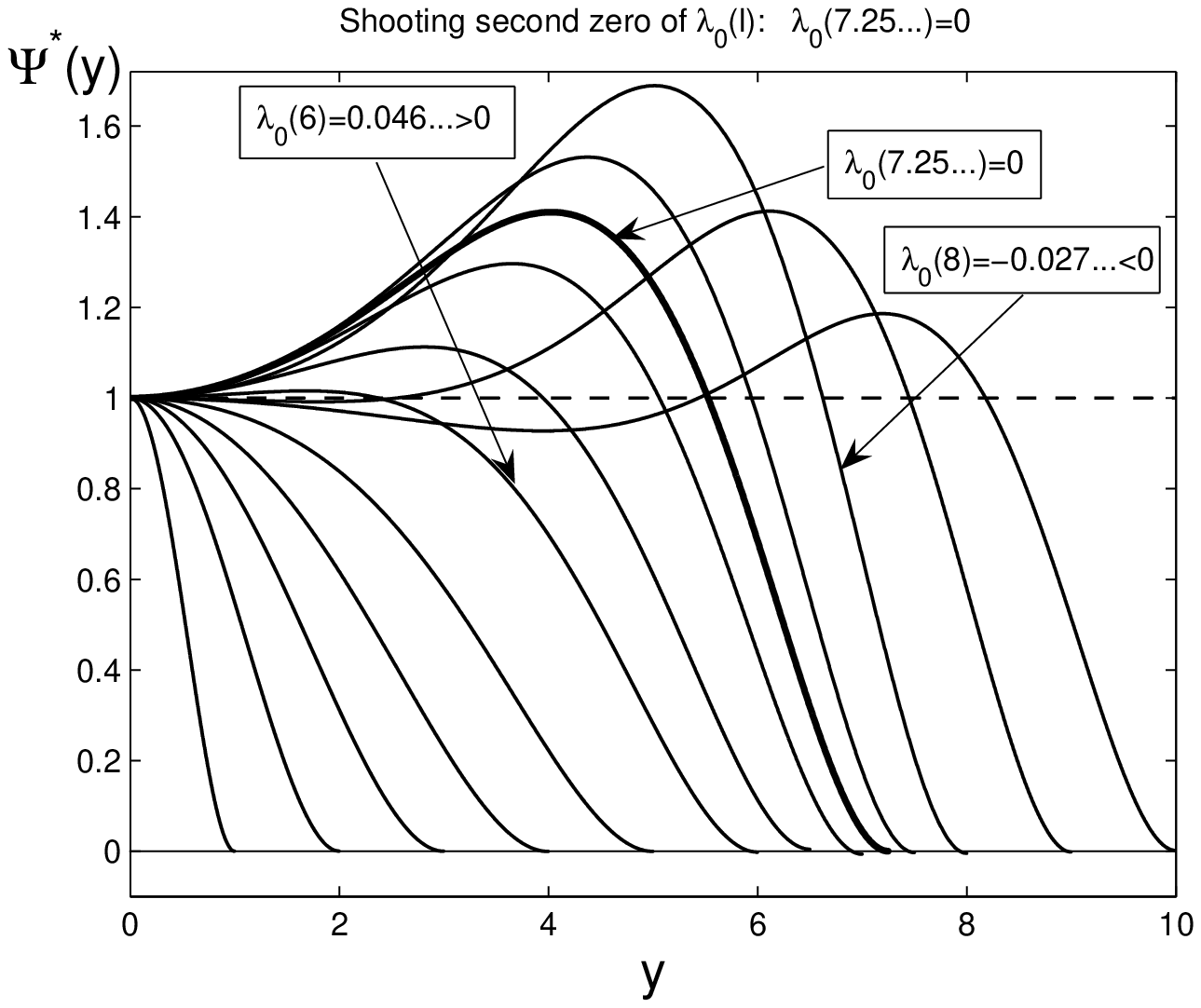}
}
 \vskip -.2cm
\caption{\rm\small Shooting zeros of $\l_0(l)$: first zero (a) and
the second one (b).}
 \label{FP2}
\end{figure}

It is clear that, for large $l$, the non-symmetric perturbation
$-\frac 14\, y D_y$ becomes essential in $\BB^*$, and Sturm's zero
property in Figure \ref{FP1} does not remain true. In Figure
\ref{FP2}(a), (b), we present the shooting results of numerical
calculating first three roots of the branch of the first
eigenvalue $\{\l_0(l)\}$:
 \be
 \label{P10}
\l_0(l)=0 \,\,\, \mbox{at first three roots} \quad l_1=4.08...\, ,
\,\,\,l_2=7.25...\, , \,\,\, \mbox{and} \,\,\,  l_3=10.\,...\, .
 \ee
 Note that $l_1$ is about $3\%$ close to the bound $l_*=3.9779...$ derived in
 \ef{P9}. Further
results of numerical experiments are shown in Table 1 fixing these
3--4 roots of $\l_0(l)$.


\begin{table}[h]
\caption{Some values of $\l_0(l)$}  
\begin{tabular}{@{}lll}
 $l$ & $\l_0(l)$ \ssk
 \\\hline
  $1$ & $-31.16...$\\ $2$ & $-1.83...$
\\ $3$ & $-0.2647...$\\
 $4$ & $-0.008152...$\\
$4.075$ & $-0.000236...$
\\
$4.0775 \approx l_1\quad$ & $0.0000113...$
\\
$4.08$ & $0.00026...$\\
 $4.1$ & $0.0022...$\\
 $4.2$ & $0.011...$\\
  $5$ & $0.0483...$\\
  $5$ & $0.0483...$\\
   $6$ & $0.046...$\\
    $7.25\approx l_2$ & $0.00167... $\\
     $7.5$ & $-0.0097...$\\
      $8$ & $-0.027...$\\
       $9$ & $-0.018...$\\
        $10\approx l_3$ & $-0.00084...$\\
         $11$ & $0.00397...$\\
          $12$ & $0.00172...$\\
           $13 \approx l_4$ & $-0.00055...$\\
\end{tabular}
\end{table}

Thus, taking into account these numerical results, in view of
\ef{P4}, we conclude that
 \be
 \label{P11}
  \fbox{$
 (0,0) \,\,\, \mbox{is regular for $l \in(0,l_1)$, is irregular
 (singular)
 for $l \in[l_1,l_2]$, etc.}
 $}
  \ee
  In addition,
since $\l_0(l_1)=0$, one can expect that for $l=l_1$, the vertex
$(0,0)$ remains irregular, but now it is not singular, i.e.,
$u(0,0^-) \not = 0$ (under a natural non-orthogonality condition)
is finite. Same happens for other roots $\l_k$ of $\l_0(l)$
(provided that the rest of eigenvalues are ``ordered", which has
not been proved for arbitrary $l>0$).

\subsection{Boundary layer
 and branching at $l=+\iy$}

\begin{figure}
\centering
\includegraphics[scale=0.8]{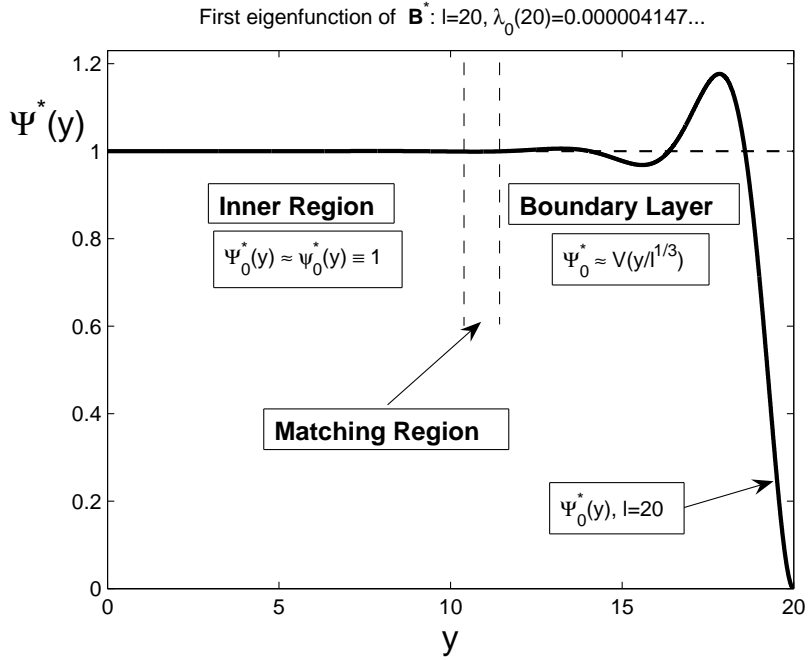} 
\vskip -.3cm \caption{\small The boundary layer structure of
$\Psi_0^*(y)$ for $l=20$.}
   \vskip -.3cm
 \label{FP3}
\end{figure}

The phenomenon of the boundary layer (BL) occurring as $l \to \iy$
is  explained by Figure \ref{FP3}. It follows that already at $l =
20$, the first eigenfunction $\Psi_0^*(y)$ of $\BB^*$ has a
 matching structure of the Boundary Layer for $y \approx l$ with the
constant eigenfunction $\psi_0^*(y) \equiv 1$ (see \ef{psi**1}),
which dominates for $l \le 12$.

The BL structure here is simple and standard; see more details in
Section \ref{S3.3}, where the BL is constructed for the parabolic
PDE. Briefly, asymptotically sharp,  in the BL for $y \sim l \gg
1$, equation \ef{P5} for $|\l_0(l)| \ll 1$ being asymptotically
small (see \ef{P17}) reads
 \be
 \label{P12}
 \begin{matrix}
 - (\Psi^*)^{(4)}- \frac 14\, l \, (\Psi^*)'+...= 0
 \LongA \Psi^*  = V(z)+...\, , \,\,\, z= l^{\frac 13} y
 \whereA\quad
 \ssk\ssk\\
 -V^{(4)}- \frac 14 \, V'=0, \,\,\, z \in (0,l^{\frac 43}); \,\,\,
 V=V'=0, \,\, z=l^{\frac 43}; \,\,\, V(0)=1.\qquad\qquad\qquad
  \end{matrix}
   \ee
The BL-function $V$ is given explicitly (cf. \ef{z12}): denoting
$b=2^{-5/3}$ and $a=\sqrt 3\, b$,
 \be
 \label{P13}
  \begin{matrix}
  V(z)=C_3 - {\mathrm e}^{-b z}\big[C_1 \cos(a z) + C_2
  \sin(a z)\big], \,\,\,
  C_1= \frac{{\mathrm e}^{b l^{4/3}}[\cos(a l^{4/3})- \frac
  1{\sqrt 3} \,\sin (a l^{4/3})]}
  {1- {\mathrm e}^{b l^{4/3}} [\cos(a l^{4/3})- \frac
  1{\sqrt 3}\, \sin (a l^{4/3}) ]}, \quad\ssk\\
 C_2= \frac{{\mathrm e}^{b l^{4/3}} [\sin(a l^{4/3})+ \frac
  1{\sqrt 3}\, \cos (a l^{4/3}) ]}
  {1- {\mathrm e}^{b l^{4/3}} [\cos(a l^{4/3})- \frac
  1{\sqrt 3}\, \sin (a l^{4/3}) ]},\quad
   C_3= \frac {1}
  {1- {\mathrm e}^{b l^{4/3}} [\cos(a l^{4/3})- \frac
  1{\sqrt 3}\, \sin (a l^{4/3}) ]}.\quad
   \end{matrix}
   \ee
The eigenvalue $\l_0(l)$ for $l \gg 1$ is then obtained by
matching in $\re$: extending $\hat \Psi_0^*= \Psi_0^* H(l-y)$,
where $H$ is the Heaviside function, we have (see Section
\ref{S3.4} for details)
  \be
  \label{P14}
  \begin{matrix}
 (\hat \Psi_0^*)^{(4)}(y)= (\Psi_0^*)^{(4)}(y) H - (\Psi_0^*)'''(l) \d(y-l) -
 (\Psi_0^*)''(l) \d'(y-l),\qquad\qquad\ssk\ssk\\
 \mbox{where} \,\,\, (\Psi_0^*)''(l)= l^{\frac 23}
 V''(l^{\frac 43})+..., \quad (\Psi_0^*)'''(l)= l \,
 V'''(l^{\frac 43})+...\,.\qquad\qquad\quad
  \end{matrix}
  \ee
 Then \ef{P5} for $\l_0(l)$ reads
  \be
  \label{P15}
   \BB^* \hat \Psi^* - l\, V'''(l^{\frac 43})\d(y-l) - l^{\frac 23}
   V''(l^{\frac 43})\d'(y-l)+... = \l_0(l) \hat \Psi_0^* \inB \re.
   \ee
Finally, substituting $\hat \Psi_0^*=\psi_0^* + w \equiv 1+ w$,
where $ w \bot
  \psi_0$, yields
 \be
 \label{P151}
 \BB^* w - l\, V'''(l^{\frac 43})\d(y-l) - l^{\frac 23}
   V''(l^{\frac 43})\d'(y-l)+...= \l_0(l)+...\, .
   \ee
 Hence, by the
orthogonality condition to $\psi_0(y)=F(y)$ (recall that $\BB^* w
\bot \psi_0$), \ef{P151} implies the following asymptotic
expression for the eigenvalue:
 \be
 \label{P16}
 \begin{matrix}
 \l_0(l)=
 -\langle l\, V'''(l^{\frac 43})\d(y-l) - l^{\frac 23}
   V''(l^{\frac 43})\d'(y-l), \psi_0 \rangle +... \ssk\ssk\\
 \equiv - l\, V'''(l^{\frac 43}) \psi_0(l) + l^{\frac 23}
   V''(l^{\frac 43}) \psi_0'(l)+... \, . \qquad\quad
    \end{matrix}
    \ee
Overall, using the asymptotics \ef{i5}, \ef{i6} for the rescaled
kernel and taking into account \ef{P13}, we obtain from \ef{P16}
the following typical approximate oscillatory behaviour:
 \be
 \label{P17}
  \tex{
  \l_0(l) \sim l^{\frac 23} {\mathrm e}^{-\hat d_0 \, l^{4/3}}
  \cos\big( \hat b_0 \, l^{\frac 43}\big) \forA l \gg 1,
  \quad \hat d_0 \sim d_0+b, \,\, \hat b_0 \sim \frac{b_0+a}2,
   }
   \ee
   where we omit constants.
 In fact, \ef{P17} shows how an $l$-branch of the first eigenvalue
 $\l_0(l)$ bifurcates at $l=+\iy$ from $\l_0(\iy)=0$ for the
 operator $\BB^*$ with the spectrum \ef{spec1}, $m=2$. It is not
 that difficult to show that similar branching occurs from any
 eigenvalue
  \be
  \label{P18}
   \tex{
   \l_k=- \frac k4 \forA k=0,1,2,... \LongA \exists \,\,\, \l_k(l)
   \forA l \gg 1, \,\,\,  \l_k(\iy)=- \frac k4,
    }
    \ee
    where \ef{P17} is evidently the only eigenvalue that can
    change sign.
The branching happens in the framework of classic perturbation
theory for linear operators (see Kato \cite{Kato}), so we do not
treat this any further. Actually, the branching \ef{P18} or
\ef{771} is the only for us reason
 for the operator \ef{P3} to admit many real
eigenvalues for various values of\footnote{We then state an open
problem: is there any direct proof that \ef{P3} has real
eigenvalues only, and the eigenfunctions form a complete and
closed set in $L^2(I_l)$? Or this is hopeless and branching
together with global extensions of the branches (also an open
problem in general) is the only reason.} $l>0$.

\section{Bi-harmonic PDE for $m=2$: Blow-up asymptotic derivation of Petrovskii-type
regularity criterion by eigenfunction expansion}
 \label{S3}

Regardless nonexistence of a definite answer concerning regularity
of the vertex of the fundamental parabolae \ef{P1}, nevertheless,
we next show that there is a way to move towards a Petrovskii-like
criterion for the case of more expanding lateral boundaries for
functions \ef{vv1}. However, the results in Section \ref{S.l}
(treated according to the limit $l=l(\t) \to +\iy$) imply, in view
of the oscillatory behaviour of the first eigenvalue in \ef{P17},
that a direct regularity conclusion for a given $\var(\t)$ is not
possible in principle. To get the regular vertex $(0,0)$, a
procedure of {\em oscillatory cut off of the boundary}
 must be performed.

At this moment, in
Appendix A,
 we continue our study with  a  certain brief, rather
questionable, and even controversial suggestion, which nevertheless
is quite attractive and actually goes along the lines of
Petrovskii's barrier ideas, so we fairly
 believe this is a right (and unique) place for such a discussion.


 As a next step, we   return
 to the spectral methods that cover both regular/irregular issues for \ef{ir12},
 which
 indeed are more complicated than for the heat equation \ef{ir1}.

\subsection{Two-region expansion}

Thus, we return to the rescaled problem \ef{ph4}  for $m=2$.
 As usual in any matching asymptotic analysis, this blow-up
 problem is solved by matching of expansions in two regions:

 \ssk

 (i) {\em Inner Region}, which is situated around the origin $y=0$,  and

 \ssk

 (ii) {\em Boundary Region} close to the boundaries $y= \pm \var(\t)$, where a boundary layer
 occurs.

 \ssk

 \noi Actually, such a two-region structure, with the
 asymptotics specified below, defines the class of generic solutions under
 consideration according to \ef{gen} (those in \ef{genNN} are all different).
We begin with the simpler analysis in the Boundary Region (ii).

 \subsection{Boundary Layer (BL) structure}
   \label{S3.3}

 Sufficiently close to the lateral boundary of $Q_0$, it is
 natural to introduces the variables
  \be
  \label{z1}
   \tex{
   z= \frac y{\var(\t)} \andA v(y,\t)=w(z,\t)
    \LongA w_\t= - \frac 1{\var^4} \, w_{zzzz} - \frac 14\, z w_z +
    \frac {\var'}\var \, z w_z.
    }
    \ee

We next introduce the BL-variables
  \be
  \label{z2}
  \tex{
  \xi= \var^{\frac 43}(\t)(1-z), \quad \var^{\frac 43}(\t) {\mathrm d} \t={\mathrm d}s,
   \andA w(z,\t)= \rho(s) g(\xi,s),
  }
  \ee
  where $\rho(s)$ is an unknown slow decaying (in the same natural sense,
 associated with \ef{vv2})
   time-factor depending on the
  function $\var(\t)$, e.g.,   $\sim \frac 1{\var(\t)}$ as a clue.
On substitution into the PDE in \ef{z1}, we obtain the following
perturbed equation:
 \be
 \label{z3}
  \begin{matrix}
  g_s= \AAA g - \frac 14\, \var^{-\frac 43} \xi g_\xi - \frac
  {\var'_\t}{\var}\, \big(1-\xi \var^{-\frac 43}\big) g_\xi \ssk\ssk\\
- \frac 43\, \var^{-\frac 13}
  \var'_\t \xi g_\xi -  \frac {\rho'_s}{\rho} \,
  g, \quad  \mbox{where} \quad  \AAA g= - g^{(4)} + \frac 14\, g'.
    \end{matrix}
   \ee
As usual in boundary layer theory\footnote{It was always  key for
PDE theory.  Shortly after the Blasius construction (1908) of the
exact
 self-similar solution for the two-dimensional  boundary layer
 equations proposed by Prandtl in 1904 (see references in \cite[p.~48]{GSVR}), similarity solutions of
  linear and nonlinear  boundary-value problems became more and more common in the
 literature.}, following  \ef{gen}, we are looking for a generic
 pattern of the behaviour described by \ef{z3} on compact subsets
 near the lateral boundary,
  \be
  \label{z4}
  |\xi| = o\big(\var^{-\frac 43}(\t)\big)
   \LongA |z-1| = o\big(\var^{-\frac 83}(\t)\big) \asA \t \to +
   \iy.
   \ee
On these space-time compact subsets, the second term on the
right-hand side of  \ef{z3} becomes asymptotically small, while
all the others are much smaller in view of the slow growth/decay
assumptions such as \ef{vv2} for $\var(\t)$ and $\rho(s)$.

\begin{figure}
\centering
\includegraphics[scale=0.65]{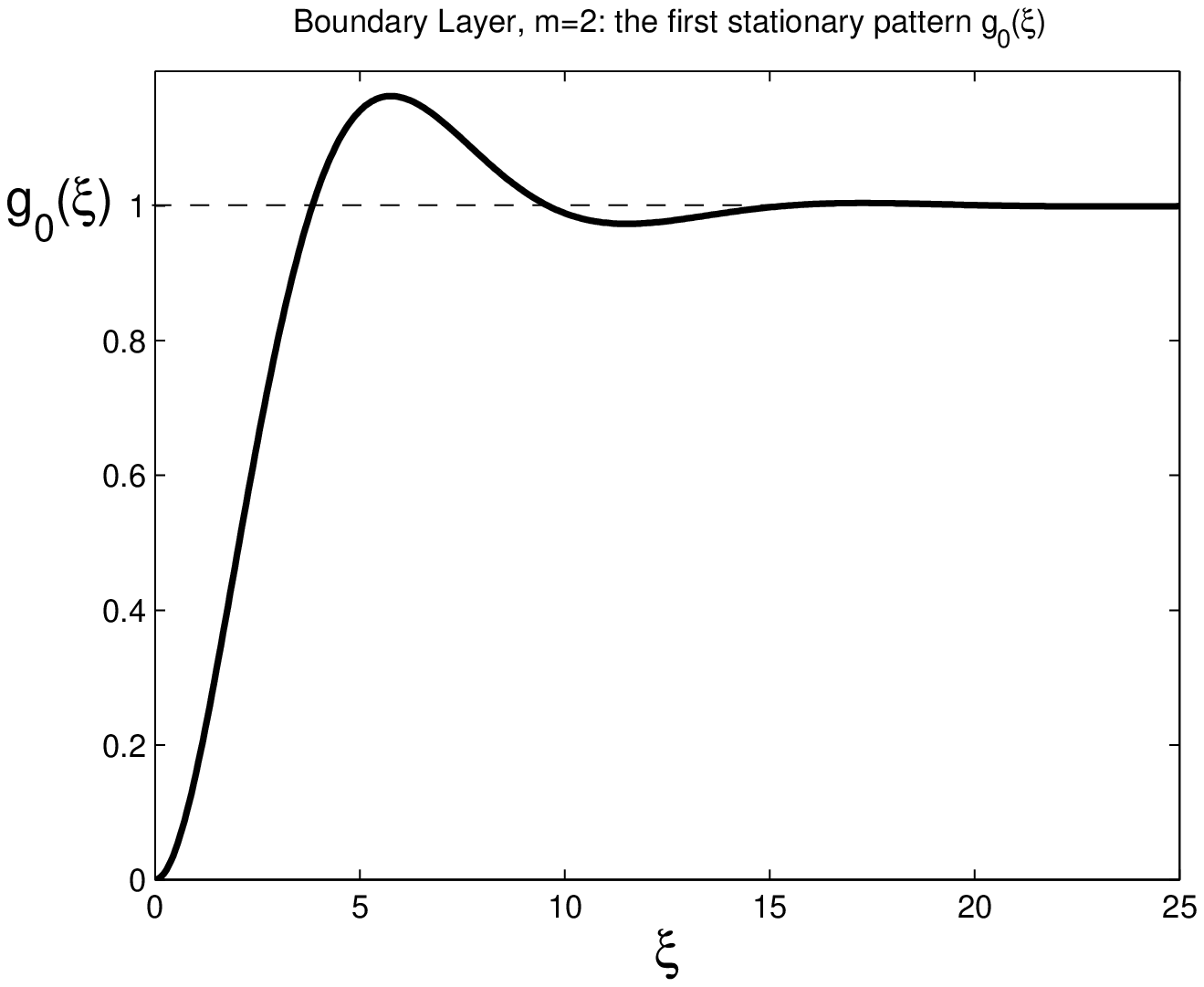} 
\vskip -.3cm \caption{\small The unique stationary solution
$g_0(\xi)$ of the problem \ef{z10}.}
   \vskip -.3cm
 \label{F1}
\end{figure}

Then posing the asymptotic behaviour at infinity:
 \be
 \label{z5}
  g(\xi,s) \to 1 \asA \xi \to + \iy \quad (\mbox{Hypothesis I for
  generic patterns}),
  \ee
  where all the derivatives are assumed to vanish, we arrive at
  the problem of passing to the limit as $s \to + \iy$ in the
  problem \ef{z3}, \ef{z5}. Assuming that, by the definition in \ef{z2},
  the rescaled orbit $\{g(s),\,\, s>0\}$ is uniformly bounded, by
  classic parabolic theory \cite{EidSys}, one can pass to the
  limit in \ef{z3} along a sequence $\{s_k\} \to +\iy$. Namely,
  by the above, we have that, uniformly on compact subsets defined in
  \ef{z4}, as $k \to \iy$,
   \be
   \label{z6}
   g(s_k+s) \to h(s) \whereA h_s=\AAA h, \quad h=h_\xi=0
   \,\,\,\mbox{at} \,\,\,\xi=0, \quad h|_{\xi=+\iy}=1.
    \ee
The {\em limit} (at $s=+\iy$) {\em equation} obtained from
\ef{z3}:
 \be
  \label{z7}
   \tex{
  h_s= \AAA h \equiv - h_{\xi\xi\xi\xi}+ \frac 14\, h_\xi
  }
   \ee
is a standard linear parabolic PDE in the unbounded domain
$\re_+$, though it is governed by a non self-adjoint operator
 $\AAA$. We need to show that, in an appropriate weighted
 $L^2$-space if necessary, under the hypothesis \ef{z5},
 the stabilization holds, i.e.,
  the $\o$-limit set of $\{h(s)\}$ consists of
 equilibria: as $s \to +\iy$,
  \be
  \label{z10}
   \left\{
   \begin{matrix}
  h(\xi,s) \to g_0(\xi) \whereA \AAA g_0=0 \,\,\, \mbox{for}
  \,\,\,\xi>0, \ssk\ssk \\
   g_0=g_0'=0 \quad \mbox{for} \quad \xi=0, \quad
  g_0(+\iy)=1.\qquad\,\,
  \end{matrix}
  \right.
   \ee
The characteristic equation for the linear operator $\AAA$ yields
 \be
 \label{z11}
  \tex{
  -\l^4 + \frac 14\, \l=0 \LongA \l_1=0 \andA \l_{2,3}= \frac 1{4^{1/3}}
  \big(-\frac 12 \pm \ii \frac{\sqrt{3}}2\big).
  }
  \ee
  This gives the unique solution of \ef{z10},  shown in Figure
  \ref{F1},
 \be
 \label{z12}
  \tex{
  g_0(\xi)=1-{\mathrm e}^{- \frac \xi{2^{5/3}}} \,\big[ \cos \big( \frac
  {\sqrt{3} \, \xi}{2^{5/3}}\big) + \frac 1{\sqrt 3} \, \sin \big( \frac
  {\sqrt{3}\, \xi}{2^{5/3}}\big) \big].
  }
  \ee

We do not concentrate on this stabilization problem \ef{z10},
which reduces to a standard spectral study of $\AAA$ in a weighted
space. Let us  mention that, in the metric of $L^2(\re_+)$, the
operator has  a clear exponential stability feature. Consider the
eigenvalue problem
 \be
 \label{bb1N}
  \tex{
   \AAA \psi \equiv  - \psi^{(4)}+ \frac 14\, \psi' =\l \psi.
    }
   \ee
As customary, multiplying  this by $\bar \psi$, next taking the
complex conjugate of the equation and multiplying by $\psi$  and
summing up the resulting equalities yields
 \be
 \label{bb2}
  \tex{
  (\l + \bar \l) \int |\psi|^2=- 2 \int |\psi''|^2 + \frac 14
  \int(\bar \psi \psi'+ \psi \bar \psi') \equiv - 2 \int
  |\psi''|^2 <0,  }
  \ee
i.e., ${\rm Re}\, \l <0$ for any $\l \in \s(\AAA)$ (this does not
take into account a continuous spectrum).

Actually, the convergence \ef{z6} and \ef{z10} for the perturbed
dynamical system \ef{z3} is the main {\sc Hypothesis (H)}, which
characterizes the class of generic patterns under consideration,
and then \ef{z5} is its partial consequence. Note that the uniform
stability of the stationary point $g_0$ in the limit autonomous
system \ef{z7} in a suitable metric will guarantee that the
asymptotically small perturbations do not affect the omega-limit
set; see \cite[Ch.~1]{AMGV}.

We stop further discussion concerning the passage to the limit in
\ef{z3} and  summarize the conclusions as follows:

\begin{proposition}
 \label{Pr.g0}


Under the given hypothesis and conditions, the
 problem $\ef{z3}$ admits a family of solutions (called generic)
 that satisfy $\ef{z10}$.

 \end{proposition}

We must admit that such a definition of generic patterns looks
rather non-constructive, which is unavoidable  for higher-order
parabolic PDEs without positivity issues. Note that,  according to
\ef{Pplus}, a similar obscure issue appears even for the heat
equation in \ef{ir1}. One can expect that \ef{z10} occurs for
``almost all" solutions, excluding just those that have eventually
the faster vanishing first Fourier coefficient that some others.
Then another BL is needed, but we do not intent to describe an
invariant manifold structure of   such a thin solution set (recall
that \ef{z3} is a difficult {\em non-autonomous} PDE).

\subsection{Inner Region analysis: towards the regularity
criterion}
 \label{S3.4}

In Inner Region, we deal with the original rescaled problem
 \ef{ph4}.
Without loss of generality, for simplicity of key calculations, we
 consider symmetric solutions defined  for $y>0$ by
 assuming  the symmetry conditions:
   \be
   \label{ss111}
  v_y=v_{yyy}=0 \quad \mbox{at\,\,\, $y=0$}.
 \ee

  In order to apply the standard eigenfunction expansion
 techniques by using the orthonormal set of polynomial
 eigenfunctions of $\BB^*$ given in \ef{psi**1}, as customary in classic PDE theory, we extend
 $v(y,\t)$ by 0 for $y > \var(\t)$:
  \be
  \label{a1}
 \hat v(y,\t) = v(y,\t)H(\var(\t)-y)=
  \left\{
   \begin{matrix}
   v(y,\t) \forA 0 \le y < \var(\t), \\
\,\, 0  \,\,\forA \,\,  y \ge  \var(\t),
\end{matrix}
 \right.
 \ee
 where $H$ is the Heaviside function.
Since $v=v_y=0$ on the lateral boundary $\{y= \var(\t)\}$, one can
check that, in the sense of distributions,
 \be
 \label{a2}
  \begin{matrix}
 \hat v_\t= v_\t H, \quad \hat v_y= v_y H, \quad
 \hat v_{yy}= v_{yy} H, \ssk\ssk\\
 \hat v_{yyy}= v_{yyy} H - v_{yy}\big|_{y=\var}\d(y-\var),\ssk\ssk
 \\
 \hat v_{yyyy}= v_{yyyy} H - v_{yyy}\big|_{y=\var}\d(y-\var)-
 v_{yy}\big|_{y=\var}\d'(y-\var).
 \end{matrix}
  \ee
Therefore, $\hat v$ satisfies the following equation:
 \be
 \label{a3}
  \hat v_\t= \BB^* \hat v - v_{yyy}\big|_{y=\var}\d(y-\var)-
 v_{yy}\big|_{y=\var}\d'(y-\var) \inB \re_+\times \re_+.
  \ee
Since, obviously, the extended solution orbit \ef{a1} is uniformly
bounded in $L^2_{\rho^*}(\re)$, we can use the converging in the
mean (and uniformly on compact subsets in $y$) the eigenfunction
expansion via the generalized Hermite polynomials \ef{psi**1}:
 \be
 \label{a4}
  \tex{
  \hat v(y,\t)= \sum\limits_{(k \ge 0)} a_k(\t) \psi_k^*(y).
   }
   \ee
   Substituting \ef{a4} into \ef{a3} and using the orthonormality
   property \ef{Ort} yields the following dynamical system for the
   expansion coefficients: for all $k=0,1,2,...\, ,$
   \be
   \label{a5}
   \tex{
   a_k'= \l_k a_k -  v_{yyy}\big|_{y=\var(\t)} \langle
   \d(y-\var(\t)), \psi_k \rangle - v_{yy}\big|_{y=\var(\t)}
   \langle \d'(y-\var), \psi_k \rangle,
    }
    \ee
    where $\l_k= -\frac k 4$ are real eigenvalues \ef{spec1}. Recall that
    $\l_k<0$ for all $k \ge 1$. More importantly, the
    corresponding eigenfunctions $\psi_k(y)$ are unbounded and not
    monotone for $k \ge 1$ according to \ef{psi44}. Therefore,
    regardless proper asymptotics given by \ef{a5}, these inner
    patterns cannot be matched with the BL-behaviour such as
    \ef{z5}, and demand other matching theory (since these are not
    generic, the latter is not developed).


     Thus, bearing again in mind \ef{gen}, one needs to
     concentrate on the ``maximal" first Fourier generic pattern associated with
 \be
 \label{a6}
 k=0: \quad \l_0=0 \andA \psi_0^*(y) \equiv 1 \quad \big(\psi_0(y)=F(y)\big).
  \ee
 Actually, this corresponds to the centre subspace behaviour for
 the equation \ef{a5}:
  \be
  \label{a7}
  \hat v(y,\t) = a_0(\t) + w^\bot(y,\t) \whereA w^\bot \in {\rm
  Span}\{\psi_k^*, \,\, k \ge 1\},
   \ee
   and $w^\bot(y,\t)$ is then negligible relative to $a_0(\t)$.
   This is another characterization of our class of generic patterns,
   Hypothesis II.
   The equation for $a_0(\t)$ then takes the form:
\be
   \label{a8}
   \tex{
   a_0'= -  v_{yyy}\big|_{y=\var(\t)}  \psi_0(\var(\t))  + v_{yy}\big|_{y=\var(\t)}
    \psi'_0(\var(\t)).
    }
    \ee

We now return to  BL theory established the boundary behaviour
\ef{z2} for $\t \gg 1$, which for convenience we state again: in
the rescaled sense, on the given compact subsets,
 \be
 \label{9}
  \tex{
  v(y,\t) = \rho(s) g_0\big( \var^{\frac 43}(\t)(1- \frac
  y{\var(\t)}) \big)+...\,.
  }
  \ee
 By the matching of both Regions, one concludes that, for such
 generic patterns,
  \be
  \label{10}
   \tex{
   \frac {a_0(\t)}{\rho(s)} \to 1 \asA \t \to \iy
   \LongA \rho(s) = a_0(\t)(1+o(1)).
   }
   \ee
 Then the convergence \ef{9}, which by regularity is also true for the
 spatial derivatives, yields, in the natural rescaled sense,
  \be
  \label{11}
   \begin{matrix}
 v_{yy}\big|_{y=\var(\t)} \to \rho(s) \var^{\frac 23}(\t) \g_1 \to a_0(\t)
 \var^{\frac 23}(\t) \g_1, \quad \g_1= g_0''(0)>0, \ssk\ssk\\
v_{yyy}\big|_{y=\var(\t)} \to -\rho(s) \var(\t) \g_2 \to - a_0(\t)
 \var(\t) \g_2, \quad \g_2= g_0'''(0).
  \end{matrix}
   \ee
Eventually, this leads to the following asymptotic ODE for the
first expansion coefficient for generic
 patterns\footnote{This result can be stated as a theorem:
{\em For the prescribed above class of generic solutions, the
first Fourier coefficient satisfies...}\,; as we have mentioned,
we avoid such non-constructive, but rigorous, ones.}:
 \be
 \label{12}
  \tex{
   \frac {a_0'}{a_0}=  G_2(\var(\t)) \equiv \g_2\var(\t) \psi_0(\var(\t))+ \g_1
   \var^{\frac 23}(\t) \psi_0'(\var(\t))+...  \forA \t \gg 1.
   }
   \ee

For instance, this gives a first easy condition of
regularity of $(0,0)$:
  for the  generic patterns,
this is the ``negative" divergence of the following integral:
 \be
 \label{13}
  \tex{
 \int\limits^{\iy} G_2(\var(\t))\, {\mathrm d} \t \quad \mbox{diverges to $-\iy$}
 \LongA
  a_0(\t) \to 0 \asA \t \to +\iy.
 }
  \ee
 In order to obtain more practical and sufficiently sharp
 conditions of regularity, we will use the expansion \ef{i5} of
 the first eigenfunction $\psi_0(y) = F(y)$, which on substitution
 into the right-hand side in \ef{12}, where both terms are
 equivalent,
 yields
  \be
  \label{14}
   \tex{
\frac {a_0'}{a_0} =  \hat \g
   \var^{\frac 23}(\t) C_3\cos\big(b_0 \var^{\frac 43}(\t)+C_4\big)
   {\mathrm e}^{-d_0 \var^{4/3}(\t)}+...  \forA \t \gg 1,
    }
    \ee
 with some $\hat \g \not = 0$ and
  constants $C_{3,4}$ depending in an obvious way on $C_{1,2}$ in \ef{i5}
and other parameters from \ef{i6}. Integrating implies that
 \be
  \label{14NN}
   \tex{
\ln|{a_0(\t)}| = \hat \g  \int\limits^\t
   \var^{\frac 23}(s) C_1 \cos\big(b_0 \var^{\frac 43}(s)+C_2\big)
   {\mathrm e}^{-d_0 \var^{4/3}(s)} \,{\mathrm d}s +...  \forA \t \gg
   1.
    }
    \ee

\noi\underline{\em Irregularity condition}. This is
straightforward: \ef{14NN} implies that the limit of $a_0(\t)$ as
$\t \to +\iy$ can be arbitrary (i.e., not necessarily zero), if
the integral converges.
 Therefore, up to the achieved
accuracy of the expansions and matching, we state the following
condition of the irregularity of $(0,0)$ (all the irrelevant
constants are omitted):
 \be
 \label{15}
  \fbox{$
  (0,0) \quad \mbox{is irregular if} \quad \int\limits^\iy
  \var^{\frac 23}(s) \cos\big(b_0 \var^{\frac 43}(s)\big) \,
   {\mathrm e}^{-d_0 \var^{4/3}(s)}\, {\mathrm d}s \quad \mbox{converges}.
   $}
   \ee
   Then, obviously, $u(0,0^-) \not = 0$ (a natural non-orthogonality condition required) is
   finite, i.e., the vertex $(0,0)$ is irregular, but
   non-singular.
   This still looks pretty similar to Petrovskii's criterion in
   \ef{RR1}.

   \ssk

   \noi\underline{\em Regularity condition: oscillatory cut-off of $\var(\t)$}.
  Thus, assume that $\var(\t)$   is
 close  to the desired (and still unknown) {\em critical} irregularity/regularity situation,
 and that we follow the first Fourier coefficient equation \ef{14}
 or \ef{14NN}, within the given accuracy. Then,
 for using \ef{14NN} for the purpose of the regularity
 analysis, we observe that the presence in the integral  the oscillatory
 factor can always
 $ \cos(\cdot)$  violate the (uniform) divergence to $-\iy$ hypothesis in
 \ef{13}. Indeed, then $\ln |a_0(\t)| \to \pm \iy$ along different
 subsequences $\{\t=\t_k^\pm\}\to +\iy$, and the limit to $+\iy$ does not
 guarantee \ef{13} (of course, for such limits, our matching with
 the BL can be violated as well, but we do not discuss this issue
 here). In other words, in the case of oscillatory kernels, a {\em
 pure divergence of the integral in \ef{14NN} does not guarantee
 the point regularity.} This is the principal difference with the
 positive kernel case; cf. \ef{993} below with no oscillatory
 component in the integral. Let us note that, in general, the
 standard divergence of \ef{14NN} will ensure that that the vertex
 $(0,0^-)$ is irregular and is also singular in the sense of
 \ef{P99}. Moreover, $u(x,t)$ has {\em infinite oscillatory} behaviour
 as $(x,t) \to (0,0^-)$, i.e., eventually oscillations get both
 $\pm \infty$.

 Thus, under the accepted hypothesis, the oscillatory part
 of the rescaled fundamental kernel \ef{i5} of the operator under
 consideration can violate the regularity of boundary point for
 {\em any} slow growing factor $\var(\t)$ satisfying
 \ef{vv1}--\ef{al1}, even if the integral diverges.
 As we have seen in Section \ref{S.l},
 for higher-order operators,
 the regularity condition ``above"
 the fundamental parabolae \ef{P1}
 becomes rather subtle and sensitive. Indeed, according to the
 integral in \ef{14NN}, arbitrarily small perturbations of
 $\var(\t)$  can ``switch over" regularity to the irregularity and {\em vice
 versa}.

 To guarantee regularity, an extra ``oscillatory cut-off" of
 $\var(\t)$ is then necessary to be introduced after a necessary
 remark.

\ssk

\subsection{
Remark 1: an analogy with  elliptic
   theory}
    This principal difficulty concerning the {\em kernels of changing sign}
    has a known
counterpart in regularity elliptic theory for \ef{LL1}. Namely, it
was first shown in \cite{MN86} (1986) for $m=2$, $N \ge 8$ that
the vertex of a cone can be irregular (singular) if the
fundamental solution of $L(\partial)$ for
 $N>2m$\footnote{For $N=2m$, the structure of the fundamental
 solution is different and more
  signed-determined:
  $$
  \tex{
F(x)= \kappa \ln \frac 1{|x|}+ \Psi\big(\frac x{|x|}\big)
\quad(\kappa(N)={\rm const}.),
 }
 $$
that makes the case special admitting arbitrary operators
$L(\partial)$ \cite[\S~8]{Maz02}; Wiener's  $m=1$, $N=2$
included.}
 \be
 \label{LL2}
  \tex{
  F(x)= \Psi\big(\frac x{|x|}\big) |x|^{2m-N}, \quad x \in \ren \setminus\{0\} \quad (N>2m)
  }
  \ee
  {\em changes sign}; see further comments and references in
  \cite[\S~1]{Maz02}. Then in elliptic theory, the regularity
  analysis is performed for a restricted subclass of operators
  $L(\partial)$ called {\em positive with weight $F$}. This actually means that $F>0$
  \cite[\S~3]{Maz02} (cf. a positivity-like condition in \cite{Freh06} for $m=2$ and \cite{Eil01}),
   so, until now,  refined elliptic regularity results
  for kernels of essentially (e.g., infinitely many) changing sign seem
  be
   unavailable.
   As in classic theory \cite{Kond67, MazPl81}, this analysis
   demands
   constructing corresponding generalized Hermite polynomials as eigenfunctions
   of an adjoint pencil of linear operators along the lines of
   those obtained for hyperbolic equations; see
   \cite[p.~254]{GSVR} and \cite{GalpSWE} (for
   $u_{tt}=-u_{xxxx}$; see the end of Section \ref{S4} below). The oscillatory cut-off of the boundary
   will then be essential. An example of such an elliptic evolution
   approach for studying singularities of \ef{DD1} in $\re^2$, $q
   \in (0,1)$
   (with many additional references presented) is given in
   \cite{BidVer99}.

   Of course,  it is well-known and this is a classic matter, that
   general theory of operator pencils
   associated with corner singularities has been well-developed for a
   number of elliptic equations;
     see \cite{KMR1, KozMaz01} and
   \cite{Freh06} for references to other papers and related monographs and as a source of further
   extensions to be traced out by the {\tt MathSciNet}.
 For other types of parabolic or elliptic problems mentioned
 above,
 new operator pencils appear.
    In the parabolic case, we do not intend to
  change the operator $-D_x^4$ in related lines\footnote{Obviously,
  unlike elliptic theory, this is not possible: no
  higher-order parabolic equations with smooth coefficients can
  have positive kernels (otherwise, comparison and the MP would be
  inherited).},
   but then need instead to
 perform a refined oscillatory cut-off of the lateral boundary.

 %


\subsection{Optimal regularity conditions for generic patterns}
\begin{figure}
\centering
\includegraphics[scale=0.75]{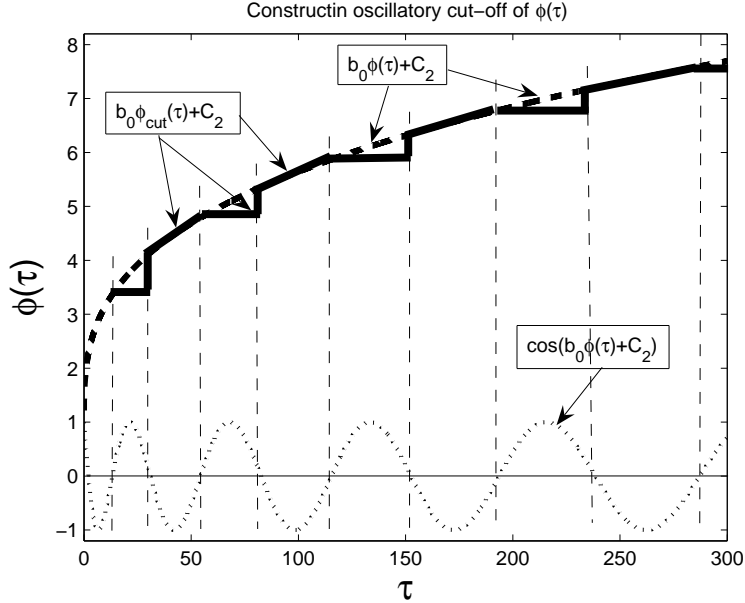} 
\vskip -.3cm \caption{\small A schematic view of a proper
oscillatory cut-off of $\var(\t)$.}
   \vskip -.3cm
 \label{FV1}
\end{figure}

  Thus,  to guarantee the regularity, we perform
   the procedure of {\em oscillatory cut-off} of a given
 function $\var(\t)$. Namely, this means that a smooth
 $\widetilde \var(\t)$ replaces $\var(\t)$, for which
  \be
  \label{RR2}
   \tex{
  \hat \g C_1 \cos\big(b_0 \widetilde{\var}^{\frac 43}(\t)+C_2\big) \approx
  \big[\hat \g C_1 \cos\big(b_0 \var^{\frac 43}(\t)+C_2\big)\big]_-
 }
   \ee
asymptotically sharp up to absolutely convergent perturbations in
the integral \ef{14NN}, where $[\cdot]_-$ denotes the negative
part. This is necessary to cut-off the positive part of the
diverging integral. Figure \ref{FV1} schematically explains how
\ef{RR2} works (for $\hat \g C_1>0$) by cutting off all positive
waves of the $\cos$-function for $\t \gg 1$ and creating ``almost
discontinuous" (jumping-like) $\widetilde \var(\t)$, which in the
figure is denoted by $\var_{\rm cut}(\t)$.
   Using necessary smoothing at the points, where $\cos(\cdot)$ changes sign,
   for any such
  monotone increasing $\var(\t)$, the corresponding $\widetilde \var(\t)$
  can be also chosen increasing. In other words, the resulting
    $\widetilde \var(\t)$ sufficiently fast jumps over
   those intervals of the length $l$, which in Section \ref{S.l} were determined
   as leading to the irregular vertex. We do not pay  special attention to the question on how such
   smoothly jumping boundary can affect the boundary layer
   behaviour described by the equation \ef{z3}, where the term
   $\var'/\var$ may then play a role (fortunately, not dominant). However, it is clear that, since
   the correct intervals of the behaviour of $\widetilde \var(\t)$ get arbitrarily long as
   $\t \to \iy$, we have enough time for establishing the
   BL-structure, so that formulae  \ef{l1} remain true closer to the end of  each of them, and
    further analysis applies. In addition, it is also clear
   that, once the solution gets very small by \ef{12} at the end
   of a fixed such interval in $\t$, a smooth monotone increasing jump of $\widetilde \var(\t)$
   at the end point cannot essentially affect the value of the solution
   during such a short interval of time. This means that, on the
   next good interval, the solution remains small and continues to be  governed
 by the eigenspace  ${\rm Span}\,\{\psi_0^*\}$.

Thus, according to the oscillatory cut-off \ef{RR2}, the condition
of regularity is as follows:
\be
 \label{15+}
  \fbox{$
  (0,0) \quad \mbox{is regular if} \quad \int\limits^\iy
 \widetilde{ \var}^{\frac 23}(s) \cos\big(b_0 \widetilde{\var}^{\frac 43}(s)+C_2\big) \,
   {\mathrm e}^{-d_0 \widetilde{\var}^{4/3}(s)}\, {\mathrm d}s \quad
   \mbox{diverges}.
   $}
   \ee
There is a clear gap between the irregularity \ef{15} and the
regularity \ef{15+} conditions, so that an Osgood--Dini-like
criterion cannot be available in principle.




 With such a procedure,
 similar to Petrovskii's presentation \ef{PP1} of the boundary regularity
results, we deduce that:
\be
 \label{16}
  \fbox{$
 \begin{matrix}
 {\rm (i)} \,\, \widetilde R(t)= 3^{-\frac 34}\, 2^{\frac {11}4}\,
 (-t)^{\frac 14}\,\big[\ln|\ln(-t)|\big]^{\frac 34}
\,\Longrightarrow\, (0,0) \quad \mbox{is regular}, \quad
\mbox{and} \qquad\quad\ssk\ssk\\ {\rm (ii)}\,\, R(t)= \big
(3^{-\frac 34}\, 2^{\frac {11}4}+\e\big)\,(-t)^{\frac
14}\,\big[\ln|\ln(-t)|\big]^{\frac 34}, \,\, \e>0
 \,\Longrightarrow\, (0,0) \quad \mbox{is irregular},
  \end{matrix}
  $}
   \ee
where, in  the regular $\widetilde R(t)$, the oscillatory cut-off
meaning \ef{RR2} must be assumed.


\subsection{Remark 2: on other asymptotic patterns and related
regularity}
 It follows from the dynamical system \ef{a5}, that, in
general, there  exist other (non-generic) asymptotic patterns
corresponding to the behaviour on each of a 1D stable eigenspace
of $\BB^*$,  which hence are not governed by the expansion
\ef{a7}. Such a behaviour then will generate its own
regularity/irregularity conditions with or without cut-offs of the
corresponding similar integrals. Indeed, those asymptotic patterns
will demand a different ``BL-like" theory, for which \ef{z5} is
not true.
 Since these are not generic and belong to the
subspace of co-dimension one, a correct posing of the IBVP in
$Q_0$ for such regular solutions is quite tricky and demands {\em
a priori} unknown and non-constructive conditions on initial data
$u_0$, so we have no reasons to take those into account.
On the other hand, describing countable sets of various blow-up
singularities is a serious  problem of modern PDE theory, but
currently this has a little to do with the boundary regularity
analysis.

\subsection{Remark 3: for the heat equation, the
spectral--BL blow-up approach is sharp}
 \label{S3.8}

Let us very briefly repeat and list the main steps of the above
blow-up analysis for the classic problem \ef{ir1} to prove
\ef{genPet}. Thus, \ef{ph1} reads
 $$
 R(t)=(-t)^{\frac 12} \var(\t),
  $$
  where, in \ef{ph4}, we get the classic Hermite operator
  \cite[p.~48]{BS}
   $$
   \tex{
    \BB^*= D_y^2- \frac 12 \, y D_y
    \inB L^2_{\rho^*}(\re), \quad \rho^*(y)={\mathrm e}^{-\frac{y^2}4},
    \quad \s(\BB^*)=\big\{ \l_k=- \frac k2, \,\,
    k=0,1,2,...\big\},
    }
    $$
    with ${\mathcal D}(\BB^*)= H^2_{\rho^*}(\re)$, etc.
Similar to \ef{z2},  the BL variables and asymptotics are now
   \be
   \label{991}
    \begin{matrix}
   \xi= \var^2(1-z) \equiv \var(\var-y) \andA
   w(z,\t)=a_0(\t)g(\xi,\t)+...\,, \qquad\quad\ssk\ssk\ssk\\
\mbox{where} \,\,  g(\xi,\t) \to g_0(\xi), \,\,\, \t \to \iy \,\,
\mbox{and} \,\,
    g_0(\xi)=1-{\mathrm e}^{-
   \xi /2} \,\,\big(g_0''+ \frac 12\, g_0'=0\big).\qquad\quad
   \end{matrix}
    \ee
 Using two formulae in the first line in \ef{a2} with
  $$
  \hat v_{yy}= v_{yy} H- v_y\big|_{y=\var} \d(y-\var),
 $$
 we arrive at the equation (the analogy of \ef{a3})
  $$
  \hat v_\t= \BB^* \hat v + v_y\big|_{y=\var} \d(y-\var).
   $$
   Hence, asymptotically via the BL structure \ef{991}, for the generic patterns as usual, the
   analogy of \ef{l2} takes much simpler form:
    \be
    \label{992}
     \tex{
 \frac{a_0'}{a_0}=-\g_1 \var \,\psi_0(\var)+...
 \LongA \frac{a_0'}{a_0}= -\frac{\g_1}{2 \sqrt \pi}\, \var\, {\mathrm e}^{-
 \frac{\var^2}4}+...\,,
  }
  \ee
  where $\g_1= g_0'(0)= \frac 12$ and $\psi_0(y)=F(y) = \frac 1{2 \sqrt \pi}\,{\mathrm e}^{-y^2/4}$
  is the positive rescaled Gaussian kernel of the heat operator. Thus, \ef{992} justifies that,
  for the generic patterns, the regularity criterion reads
   \be
   \label{993}
    \fbox{$
\mbox{$(0,0)$ is regular iff} \quad
 \int\limits^\iy \var(\t)\, {\mathrm e}^{-
 \frac{\var^2(\t)}4}\, {\mathrm d}\t = +\iy
  $}
  \ee
  (no cut-off is necessary in this non-oscillatory case).
One can see that \ef{993} is equivalent to that in \ef{RR1} or
\ef{999}. Recall that there exist other asymptotic patterns
related to the stable subspace of $\BB^*$ of co-dimension 1,
which, as usual, is not taken into account.

\section{On extensions to $2m$th-order poly-harmonic and other
operators}
 \label{S4}

We discuss the extensions of the asymptotic method to other PDEs.

\subsection{Poly-harmonic equation}

For \ef{PE1}, the  spectral properties of the operator
 \be
 \label{41}
  \tex{
 \BB^*= -(-\D)^m - \frac 1{2m} \, y \cdot \n
 }
  \ee
 are given in \cite{Eg4}. A sharp asymptotic expansion (similar to \ef{i5}) of the radially
 symmetric rescaled kernel $F(|y|)$ satisfying
 \be
  \label{m1}
  \tex{
 \BB F=0 \inB \ren, \quad \int F=1 \LongA
 \psi_\b(y)= \frac{(-1)^{|\b|}}{\sqrt{\b !}}\, D^\b F(y), \,\,\,
 |\b|=0,1,2,...\, ,
 }
 \ee
 where $\b$ is a multiindex in $\ren$,
  is also not that difficult to
 get. Then \ef{m1} implies that the adjoint basis $\{\psi_\b^*\}$ consists
 of generalized Hermite polynomials.

 In the radial case,  the critical $R(t)$ for boundary regularity is about
  \be
  \label{42}
   R(t)= C_*\, (-t)^{\frac 1{2m}}\, \big[\,\ln |\ln (-t)|\,\big]^{\frac{2m-1}{2m}} \whereA
   C_*=d_0^{-\frac 1 \a} >0 \,\,\, \mbox{is a constant},
    \ee
and the oscillatory cut-off is assumed for the regularity. As
usual,  replacing $C_*$ by $C_*+\e$ for any $\e>0$ makes $(0,0)$
irregular.

  Some calculations  become more involved, especially
 in the case of non-radial $Q_0$, where the spatial shape of the
 shrinking cusp will affect the right-hand side of
  \ef{14}, where the BL structure, though having similar variables such as
  \ef{z2}, is now governed by complicated elliptic problems
  instead  of the ODE \ef{bb1N}.
  For instance, if the boundary of $Q_0$ is given by
 the characteristic fundamental paraboloid in the $x$ and $y$
 variables, respectively (i.e., $\varphi(\t) \equiv 1$):
  \be
  \label{par1}
   \tex{
  \sum_{(|\a|=2m)} a_\a x^\a=(-t) \LongA
  \sum_{(|\a|=2m)} a_\a y^\a=1.
  }
  \ee
The case  $ a_\a>0$ for any $\a$ corresponds to Miha${\rm
\check{i}}$lov's  \cite{Mih61, Mih63I, Mih63II} and Kondrat'ev's
  \cite{Kond66} cases.
Hence, the regularity in the strong sense of the vertex $(0,0)$
 depends on the spectrum of the operator \ef{41} in the domain
 $I_a$
 with the boundary given by \ef{par1}.
   Similar to the analysis in Section
  \ref{S.l}, we have that, by Poincare's inequality, $(0,0)$ is
  regular if  the diameter of $I_a$ is not that large
  (since the eigenvalues of $\BB^*$ are close to those for $-(-\D)^m <0$). On the contrary,
  since $\BB^*$ is not self-adjoint and sign-definite, {\em there are coefficients
  $\{a_\a>0\}$, for which the vertex $(0,0)$ is not regular in the classic sense.}

If some of  the coefficients $\{a_\a\}$ are negative (this is
 Fe${\rm \check{i}}$gin's case \cite{Fei71}),  $\BB^*$ is then
 posed in an unbounded domain with a non-compact boundary, and its
 spectral theory becomes more involved. Nevertheless, for such
 sufficiently ``thin" domains, we may expect the vertex be regular,
 and sometimes
 irregular otherwise.



\subsection{Linear dispersion equations}

There is no much  difference in studying the boundary regularity
for odd-order linear PDEs. For instance, as a simple example,
consider the third-order {\em linear dispersion equation} in a
similar setting:
 \be
 \label{l1}
   u_t=u_{xxx} \inB Q_0=\{-R_-(t) < x < R_+(t), \quad -1
 <t <0\},
  \ee
  where $R_\pm(t)$ are positive smooth functions on $[-1,0)$,
  continuous on $[-1,0]$, and $R_\pm(0)=0$.
  Since the operator $D_x^3$ is anisotropic, the functions
  $R_\pm(t)$ are essentially different.
 According to odd-order PDE theory (see e.g., \cite{Fam02}), for \ef{l1}, it is
 allowed to put the following Dirichlet boundary conditions:
 \be
 \label{l2}
 u=u_x=0 \atA x=R_-(t) \andA u=0 \atA x=R_+(t), \quad t \in
 [-1,0).
  \ee
 Since, as we will show, the main singular phenomena occur at the
 right-hand ``oscillatory" lateral boundary $\{x=R_+(t)\}$, we
 concentrate on   this analysis, and, in general, can put
  \be
  \label{R-}
  R_-(t) \equiv - \iy,
  \ee
  so excluding this part of the boundary from consideration. Then,
the regularity criterion of $(0,0)$ will solely depend on the
behaviour of $R_+(t)$ as $t \to 0^-$.

Thus, taking into account the right-hand lateral boundary, as in
\ef{ph1}, we perform the first rescaling:
 \be
 \label{ph4N}
  \begin{matrix}
 R_+(t)=(-t)^{\frac 13} \var_+(\t), \quad y= \frac
 x{(-t)^{1/3}}, \quad \t=-\ln(-t), \quad u(x,t)=v(y,\t),
 \ssk\ssk\\
  \mbox{where} \quad
  v_\t= \BB^* v \equiv  v_{yyy}- \frac 13 \, y v_y, \quad
  v=0 \atA y= \var_+(\t), \,\,\t \ge 0.\quad\,\,
   \end{matrix}
  \ee
  Hence,
unlike parabolic equations,
 for (\ref{l1}),  another (but indeed similar in many places)
version of Hermitian spectral theory for $\BB^*$ is necessary to
get optimal
 behaviour of $R_+(t)$ as $t \to 0^-$ for the
 regularity/irregularity of $(0,0)$.

First, as above, the regularity criterion is associated with the
structure of
 the fundamental solution of \ef{l1} that has the similarity form
  \be
  \label{gg1}
   \tex{
  b(x,t)=t^{-\frac 13}F(y),
  \quad y=\frac  x{t^{1/3}},
  }
  \ee
 where  $F={\rm Ai}(y)$
 is the  Airy function satisfying
    \be
    \label{ai1}
     \mbox{$
 {\bf B}F \equiv   F''' + \frac 13 \, (y F)'=0 \quad \mbox{in} \quad \re, \quad \int
    F=1 \LongA F'' + \frac 13\, y F=0.
    $}
     \ee
 Hence, the rescaled kernel $F(y)$ has an exponential decay as
 $y \to - \infty$ only, and as $y \to +\infty$ is oscillatory
 according to the asymptotics
 \be
   \label{mmm.9}
   \mbox{$
   F(y) = \left\{
   \begin{matrix}
    C_1|y|^{-\frac 14} {\mathrm e}^{-d_0|y|^{3/2}}+... \quad \mbox{as} \,\,\, y \to - \infty,
    \qquad\qquad\qquad \quad\,\,\ssk\ssk \\
  y^{-\frac 14}\big[C_2 \cos\bigl(d_0 y^{\frac 32}\bigr
) + C_3 \sin\bigl(d_0 y^{\frac 32}\bigr )\big]+...\quad
   \mbox{as} \,\,\, y \to + \infty,
 \end{matrix}
  \right.
 $}
  \ee
  where $d_0= \frac {2\sqrt 3}9$ and $C_{1,2,3}$ are some constants of non-essential values.

 Secondly, this (see extra details in \cite[\S~9.2]{2mSturm} and \cite{RayGI})
 generates very oscillatory as $y \to +\infty$
 and unbounded
  (for any $k \ge 1$) eigenfunctions of $\BB$ in \ef{ai1} (cf. \ef{eigen})
 \be
    \label{psiAi1}
     \mbox{$
   \psi_k(y) = \frac {(-1)^k}{\sqrt{k !}} \, D_y^k F(y) \quad
     \big(F={\rm Ai}(y)\big) \withA \s(\BB)=\big\{- \frac k3, \,\,\,
 k=0,1,2,...\big\}.
    $}
    \ee
    For existence of such a discrete spectrum, a proper ``radiation" condition
    at $y=+\iy$ is posed, which excludes  ``non-oscillatory"
    (polynomial-rational) bundles
     that are different from those in \ef{mmm.9}.
 The ``adjoint"  linear operator
 \be
 \label{bb1}
  \tex{
 \BB^*=D_y^3- \frac 13\, y D_y \withA \s(\BB^*)=\s(\BB)=\big\{- \frac k3, \,\,\,
 k=0,1,2,...\big\},
 }
 \ee
 has a complete set of eigenfunctions $\{\psi_k^*\}$ that are generalized Hermite
 polynomials. In general, then the bi-orthonormality property \ef{Ort}
 demands  using
  Hahn--Banach Theorem
 on
 extensions of  continuous linear functionals (or similar techniques of regularization of
  oscillatory integrals in
 a  v.p. or c.r. sense) and other constructions
  (this is
 not that
  crucial for in what follows).
  Both operators are defined in special weighted $L^2$
 spaces.
 One can see that
  ${\bf B}^*$ {\em is not adjoint} to ${\bf B}$ in the standard metric of $L^2$,
  since,
  obviously, then
   $\tilde {\bf B}^*=-{\bf B}+ \frac 13\, I$, so
  ${\bf B} - \frac 16 \, I$ is skew-symmetric.
 In fact,   $\BB^*$ is adjoint to $\BB$
  in a space with
 the {\em indefinite metric} given by
  \be
   \label{nn1}
  \mbox{$
 \langle v,w \rangle_*= \int v(y) \overline{w(-y)} \, {\mathrm
 d}y.
 $}
   \ee
 This case corresponds to the decomposable space with indefinite metric with straightforward
   majorizing one and is treated as rather trivial; see Azizov--Iokhvidov \cite{AI89} for  linear
  operators theory in spaces with indefinite metric (the theory was initiated in the 1940s
  and 50s by  Pontryagin
  and Krein).
Then the domain of $ \BB^*$ is defined as $H^3_{ \rho^*}$, etc.

Thus,  it is key that the bounded operator  ${\bf B}^*:
H^3_{\rho^*} \to L^2_{\rho^*}$
 admits  a complete set of  {\em polynomial} eigenfunctions
 $\Phi^*=\{\psi^*_k(y)\}$,
 which  are constructed similarly to those obtained in
 (\ref{psi**1}).
 Quite analogously, these polynomials
  are used in our regularity analysis.

\ssk

\noi\underline{\em Boundary Layer at $y=R_+(t)$}. Similar to
\ef{z1}, \ef{z2}, the BL-variables are
 \be
 \label{s1}
 \tex{
 z = \frac y{\var_+(\t)}, \quad \xi=\var_+^{\frac 32}(\t)(1-z),
 \quad w(z,\t)=\rho(\t) g(\xi,\t).
 }
 \ee
 This yields the asymptotic problem of stabilization to a unique
 profile $g_0(\xi)$ (as usual, we omit this not that easy
 analysis),
 where, as seen from \ef{ph4N},
  \be
  \label{s2}
   \tex{
    g_0''' - \frac 13\, g_0'=0, \,\,\, g_0(0)=1, \,\, g_0(+\iy)=1
    \LongA
  g_0(\xi)=1-{\mathrm e}^{-\xi/\sqrt 3}.
  }
  \ee

\ssk

\noi\underline{\em Inner Region: the regularity criterion for
$R_+(t)$}. Using the formulae as in  \ef{a2}, we obtain here the
equation
 \be
 \label{a3NN}
  \hat v_\t= \BB^* \hat v + v_{yy}\big|_{y=\var_+}\d(y-\var_+)+
 v_{y}\big|_{y=\var_+}\d'(y-\var_+) \inB \re\times \re_+,
  \ee
  and use the eigenfunction expansion \ef{a4} via the corresponding  generalized Hermite
  polynomials  that leads to the dynamical system such as \ef{a5}.
Looking for the generic patterns with the representation \ef{a7},
\ef{a3NN} yields
\be
   \label{a8N}
   \tex{
   a_0'=   v_{yy}\big|_{y=\var_+(\t)}  \psi_0(\var_+(\t))  - v_{y}\big|_{y=\var_+(\t)}
    \psi_0'(\var_+(\t)) \quad \big(\psi_0(y) \equiv F(y)\big).
    }
    \ee
By  BL theory and matching, we can use the fact that, in the
rescaled sense, on the given compact subsets, with first
derivatives
 \be
 \label{9N}
  \tex{
  v(y,\t) = a_0(\t) g_0\big( \var_+^{\frac 32}(\t)(1- \frac
  y{\var_+(\t)}) \big)+...\, .
  }
  \ee

Eventually, this leads to the following asymptotic ODE for the
first expansion coefficient for generic patterns:
 \be
 \label{12N}
  \tex{
   \frac {a_0'}{a_0}=  G_3(\var_+(\t)) \equiv \g_2\var_+(\t) \psi_0(\var_+(\t))+ \g_1
   \var_+^{\frac 12}(\t) \psi_0'(\var_+(\t))+...  \forA \t \gg 1.
   }
   \ee
 Since by the asymptotics \ef{mmm.9} both terms are similar,
 this yields
  \be
  \label{14N}
   \tex{
\frac {a_0'}{a_0} = \hat \g \,
   \var_+^{\frac 34}(\t) C_4 \cos\big(d_0 \var_+^{\frac 32}(\t)+C_5\big)
  +...  \forA \t \gg 1.
    }
    \ee

This leads to the following condition  of  irregularity of
$(0,0)$:
 \be
 \label{15N}
  \fbox{$
  (0,0) \quad \mbox{is irregular if} \quad \int\limits^\iy
  \var_+^{\frac 34}(s) \cos\big(d_0 \var_+^{\frac 32}(s)\big)
  \, {\mathrm d}s \quad \mbox{converges}.
   $}
   \ee
There is no $\log$-$\log$ factor in the critical $R_+(t)$ and just
a single $\log$
occurs: 
in view of the influence of the oscillatory factor $\cos(\cdot)$,
the point $(0,0)$
 \be
 \label{15N1}
 \fbox{$
 \mbox{for $R_+(t)=(-t)^{\frac 13}|\ln(-t)|^\g$:
 regular for $\g \le \frac 43$ and is irregular for $\g> \frac 43$,}
  $}
  \ee
 where, for the regular case, the corresponding oscillatory cut-off
 on $\var(\t)$ according to the rule as in \ef{RR2} is assumed, i.e., then
 $R_+(t)$ reads $\widetilde R_+(t)$.

\ssk

Concerning the \underline{\em left-hand lateral boundary}
$\{x=-R_-(t), \, t \in [-1,0)\}$ (setting $R_+(t) \equiv +\iy$),
the standard exponential decay behaviour of the kernel \ef{mmm.9}
as $y \to - \iy$ is supposed to give a result which is similar to
the heat equation Gaussian, with $F(y) \sim {\mathrm e}^{-y^2/4}$.
Since the decay rate is slower in \ef{mmm.9} governed by the power
$|y|^{\frac 32}$, bearing in mind the integral criterion in
\ef{15}, one can derive the following conclusion on the
regularity:
 \be
 \label{ct1}
  \tex{
  \int^\iy \,\, ...\,\,\, {\mathrm e}^{-d_0|\var_-(s)|^{3/2}}\, {\mathrm
  d}s = \iy \LongA R_-(t)=\big(\frac{3\sqrt{3}}2\big)^{\frac
  23}(-t)^{\frac 13}\big[ \ln|\ln(-t)|\big]^{\frac 23},
  }
  \ee
 where the constant $\big(\frac{3\sqrt{3}}2\big)^{
  2/3}$ cannot be increased without
 losing the regularity of $(0,0)$.
 Since the kernel $F(y)$ is of constant sign in this limit, no
 oscillatory cut-off is necessary. Thus, the third-order PDE
 \ef{l1} is the {\em last} one, for which for a (one-sided)
 regularity criterion, there occurs a non-oscillatory kernel
 and the result looks similar to the heat equation one \ef{PP1}
 (though no sub-solutions like \ef{sub1} can be used, and, hypothetically,
  a majorizing order-preserving flow
 as in
 Appendix A
 can be discussed).  As we have
 seen, for the bi-harmonic equation \ef{ir12}, both lateral
 boundaries, in general, are assumed to undergo oscillatory
 cut-offs for the regularity. The same is true for the fifth-order
 linear dispersion equation and the  tri-harmonic one
  \be
  \label{sss1}
  u_t= u_{xxxxx} \andA u_t=u_{xxxxxx}, \quad \mbox{etc.};
   \ee
 on their and others oscillatory properties, see  \cite[Ch.~4 and 3]{GSVR}.

\ssk

More complicated spectral theory of multi-dimensional odd-order
operators $\BB^*$ and $\BB$
 is
necessary for the study of the boundary regularity phenomena for
the  $N$-dimensional counterparts of linear dispersion PDEs such
as
  \be
  \label{DN1}
 u_t = (\D u)_{x_1} +... \,\,\,\,
 \mbox{or} \,\,\,\,
  u_t = -(\D^2 u)_{x_1} +...
  \quad \mbox{in} \quad  Q_0.
  \ee




\subsection{On a quasilinear fourth-order diffusion equation}

Here we briefly explain necessary changes towards the regularity
analysis for quasilinear diffusion-like PDEs. We take the 1D cubic
4th-order porous medium equation (the PME--4), as a basic model,
 \be
 \label{u1}
 u_t=-(u^3)_{xxxx} \inB Q_0,
  \ee
  where $Q_0$ is the same as in \ef{ir12} and the Dirichlet
  boundary conditions now read
 \be
 \label{u2}
 u=(u^3)_x=0 \atA x= \pm R(t), \quad -1 \le t <0.
  \ee
  We specially have fixed a fully divergent diffusion model
  \ef{u1} with the monotone operator in the metric of $H^{-2}$, for which
   existence and uniqueness is guaranteed by classic theory; see
   Lions
  \cite[Ch.~1,\,2]{LIO}.
 Using the scaling \ef{ph3} yields
 \be
 \label{u3}
  \tex{
   v_\t=\BB^*(v) \equiv -(v^3)_{yyyy} - \frac 14 \, y v_y,
    }
    \ee
    where $\BB^*(v)$ is now nonlinear. After the BL-analysis, we
    return to \ef{u3} to see how the linear operator $\BB^*$ as in
    \ef{ph4} will be recovered.

    \ssk

    \noi\underline{\em Boundary layer}. The approach is similar
    with the BL variables as in \ef{z2}, where
     \be
     \label{u4}
     \tex{
     \xi= \frac{\var^{4/3}(\t)}{\rho^{2/3}(\t)} \, (1-z),
     \quad w(z,\t)=\rho(\t) g(\xi,\t) \quad \big(z= \frac y{\var(\t)}\big).
     }
     \ee
     The BL
    rescaled equation such as \ef{z3} is easily recovered with
    the nonlinear operator
     \be
     \label{u5}
      \tex{
      g_0: \quad \AAA(g) \equiv -(g^3)^{(4)} + \frac 14\, g'=0,
      \,\,\,
       \xi>0; \quad g=(g^3)'=0\,\, \mbox{at}\,\, \xi=0, \,\,\, g(+\iy)=1.
      }
      \ee
 The proof of existence  for the problem \ef{u5} is not
 easy, but indeed doable by a standard shooting approach, though the uniqueness
  can cause a serious technical problem. To avoid an essential  declining to
   classic ODE theory, which is not the subject of the present PDE
   research, we avoid spending a couple, or possibly a few of next pages with various
   elementary and not that easy related ODE calculus.

\begin{figure}
\centering
\includegraphics[scale=0.7]{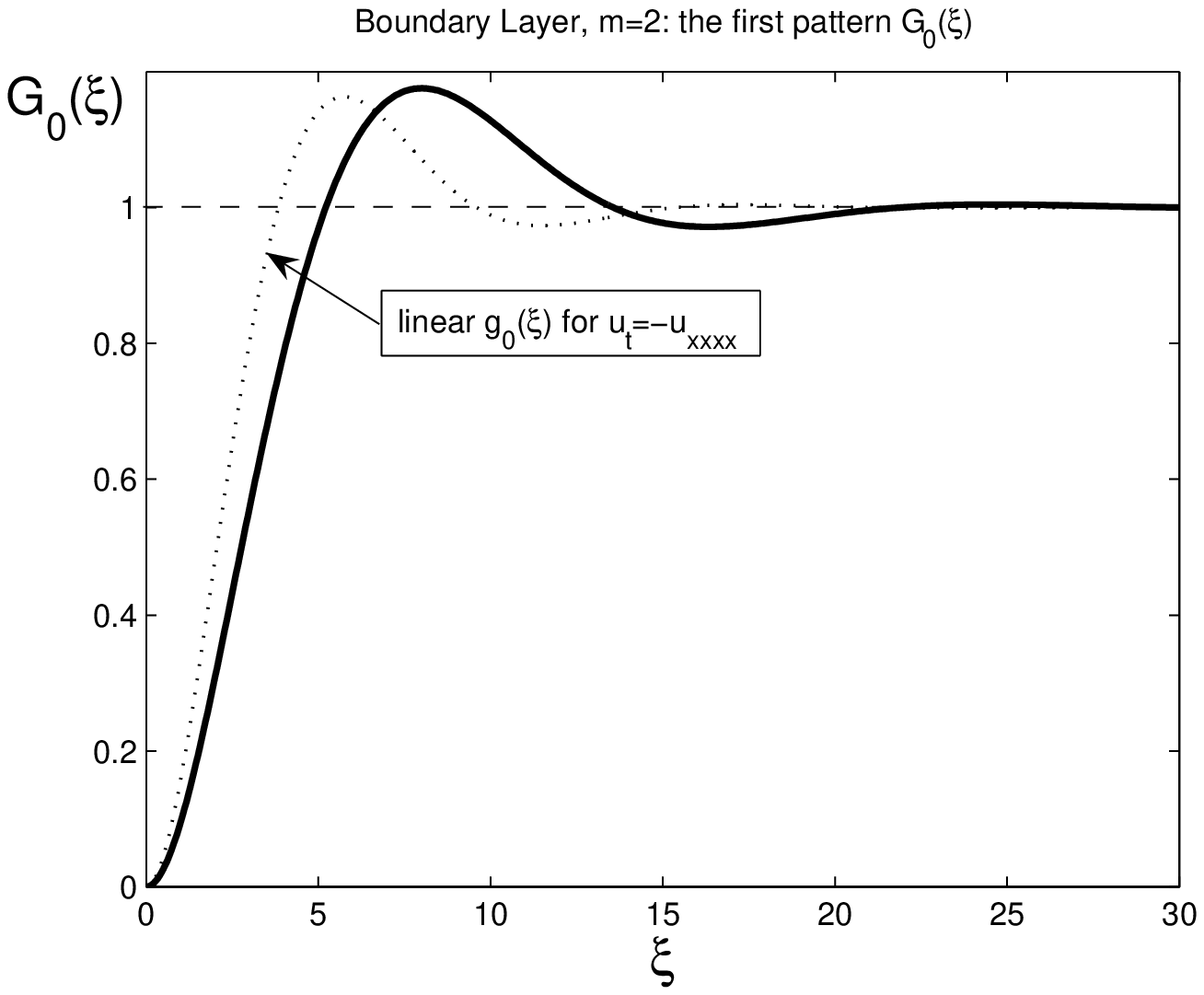} 
\vskip -.3cm \caption{\small The unique stationary solution
$G_0(\xi)$ of the problem \ef{u6}.}
   \vskip -.3cm
 \label{F2}
\end{figure}

Instead, using standard facilities,
 we again
   just check  these key properties of existence-uniqueness of $g_0(\xi)$ by
   using the solver {\tt bvp4c} of the
 {\tt MatLab}  with the enhanced  accuracy
 and tolerances \ef{Tol1}.
 The sufficient for us  and reliable numerical
 justification is shown in Figure \ref{F2}, where
  the semilinear equation is solved:
  \be
  \label{u6}
   \tex{
   G_0=g^3_0: \quad - G^{(4)} + \frac 1{12} \, |G|^{-\frac
   23} G'=0,
    \,\, G=G'=0\,\, \mbox{at} \,\, \xi=0, \,\,\, G(+\iy)=1.
      }
      \ee
      The dotted line shows the ``linear"  $g_0(\xi)$ from Figure
      \ref{F1}, so that both linear and nonlinear profiles are  similar,
      though correspond to  different ODE mathematics.

 In general, stabilization in
the full model, which is a quasilinear extension of \ef{z3}, is
rather technical and involved, but a class of such generic
solutions forming the necessary boundary layer is indeed
obtainable. Moreover, instead of the corresponding Proposition
 \ref{Pr.g0}, which is difficult to prove in a sufficient
 generality, we may view this result as a definition of the class
 of generic solutions under consideration.
 For this class of solutions, the following holds on
appropriate compact subsets near the boundary (cf. \ef{9}):
 \be
 \label{u7}
  \tex{
  v(y,\t) = \rho(\t)g_0
  \big(\frac{\var^{4/3}(\t)}{\rho^{2/3}(\t)} \, (1- \frac
  y{\var(\t)})\big)+... \asA \t \to +\iy.
   }
   \ee

 Thus, in the nonlinear case, the BL analysis is rather similar
 and demands only the extra scaling factor in the BL-variable $\xi$ in \ef{u4}; cf.
 the linear one \ef{z2}.

\ssk

\noi\underline{\em Inner Region}. This is also  similar but there
is a single nonlinear aspect that we have to address to.
 Thus, the only principal  difference between \ef{u3} and \ef{ph4} is that
  $\BB^*(v)$ is a nonlinear operator, so we cannot immediately use
  the eigenfunctions expansion \ef{a4}.

To this end, we recall that the BL-behaviour \ef{u7} actually
implies that, in this critical case, we look for the so-called
``mesa-like" asymptotics of solutions, which are supposed to be
almost flat and to be independent of $y$ for $\t \gg 1$ on compact
subsets bounded away from the BL at $|y| \sim \var(\t)$ (where
\ef{u7} takes place). Such mesa-like asymptotics
correspond to truly critical
 cases and are well-known in reaction-diffusion theory; see
 examples in \cite[p.~59,\,184]{AMGV}.
According to this, it is correct to use the following
approximation in the nonlinear term in $\BB^*( v)$:
 \be
 \label{WW1}
  v^3= \rho^2(\t) v(1+o(1)) \forA \t \gg 1
  \ee
  on compact subsets bounded away from thin domains, where the
  ``non-flat"
  BL-expansion \ef{u7} is in charge. Then
   $$
   -( v^3)_{yyyy} \mapsto - \rho^2(\t)  v_{yyyy} +...\, ,
 $$
 so that, to get back to the {\em linear operator} $\BB^*$ in the
 equation, one needs to rescale the spatial variable as follows:
   \be
   \label{lin1}
    \tex{
   v(y,\t)=  w(\eta,\t) \whereA \eta= \frac
  y{\sqrt{\rho(\t)}} < \hat \var(\t)=
  \frac{\var(\t)}{\sqrt{\rho(\t)}} \to +\iy.
  }
   \ee
 Then  $w(\eta,\t)$ solves the  equation that in leading terms
 coincides with \ef{a3},
 \be
 \label{lin4}
  \tex{
  w_\t=  \BB^* w +  \frac {\rho'}\rho\,  \frac 12\,
  \eta w_\eta +... \, ,
  }
  \ee
where $+...$ denotes higher-order terms due to the sharp
approximation \ef{WW1}. In comparison with the linear case
\ef{a3}, we observe the additional second term of the order
 $O\big(\frac {\rho'}{\rho}\big)$  for $\t \gg1$, which is clearly negligible
in view of \ef{vv1}--\ef{ph2}. Moreover, $w_\eta \equiv 0$ on the
subspace spanned by $\psi_0^*=1$, for which \ef{WW1} and \ef{z5}
hold.

We next introduce $\hat w$ according to the cut-off as in \ef{a1}
at $\eta = \hat \var(\t) \to + \iy$ given in \ef{lin1}. Then
performing it first for the quasilinear equation for $w$ and next
using the linearization as in \ef{lin4}, we arrive at
\be
 \label{lin4LL}
  \begin{matrix}
  \hat w_\t=  \BB^* \hat w +  \frac {\rho'}\rho\,  \frac 12\,
  \eta  \hat w_\eta
   \qquad
 \ssk\ssk\ssk\\
   +
   \frac 1{\rho^2}\big[
   (w^3)_{\eta\eta\eta}\big|_{\eta=\hat \var}\d(\eta-\hat
   \var)+
 (w^3)_{\eta\eta}\big|_{\eta=\hat \var}\d'(\eta-\hat \var)
 \big]
 + ... \, .
 \end{matrix}
  \ee
Using in the last terms the BL estimates following from \ef{u7}
(similar to \ef{11}), and, as usual, applying \ef{10} for the
class of generic solutions, the eigenfunction expansion \ef{a4}
for \ef{lin4LL} yields the following asymptotic ODE for the first
Fourier coefficient ($a_0>0$):
 \be
 \label{m1AA}
 \tex{
  a_0'(\t)=  \g_2 \sqrt{a_0(\t)} \,\var(\t) \psi_0\big( \frac
  {\var(\t)}{\sqrt{a_0}}\big)+ \g_1 a_0^{\frac 23}(\t)\var^{\frac 23}(\t) \psi_0'\big( \frac
  {\var(\t)}{\sqrt{a_0}}\big)+... \, .
   }
   \ee

   It follow that, unlike a simpler equation \ef{12}, due to the
   extra ``nonlinear" scaling \ef{lin1}, the Fourier coefficient
   $a_0(\t)$ in \ef{m1AA} is now entering the oscillatory rescaled
   kernel $\psi_0$ and its derivative $\psi_0'$.
 Since $\var(\t) \to +\iy$ and $a_0(\t) \to 0$ for the regularity
 (see \ef{13}), the first term in \ef{m1AA} is leading that yields
 on substitution the asymptotics \ef{i5}, \ef{i6}:
  \be
  \label{m2}
  \tex{
  a_0' = \hat \g \, a_0^{\frac 23} \var^{\frac 23}{\mathrm e}^{-
  d_0(\var/\sqrt{a_0})^{4/3}} C_3 \cos\big[ b_0 \big( \frac
  \var{\sqrt{a_0}}\big)^{\frac 43}+C_4 \big]+...\, .
   }
   \ee
The asymptotic analysis of this ODE is similar. For instance, the
following rate of convergence is obtained for the unperturbed
``Gaussian" parabola:
 \be
 \label{m4}
 \tex{
  \var(\t) \equiv 1: \quad a_0(\t) \sim 3^{\frac 32}\,
  2^{-\frac{11}2} (\ln \t)^{-\frac 32} \to 0,
   }
   \ee
which shows the obvious result that the fundamental kernel
parabolic shape for the shrinking point of $Q_0$ is regular.

In general, taking into account the main exponential term in
\ef{m2} only, i.e., fixing the simplest approximating model
 \be
 \label{m8}
  \fbox{$
a_0' = -{\mathrm e}^{-
  d_0(\var/\sqrt{a_0})^{4/3}} \quad (a_0(\t)>0),
 $}
   \ee
 we see that the critical $R(t)$ corresponds to
 \be
 \label{m5}
 \tex{
 \tilde R_*(t)
 \sim C(-t)^{\frac 14}\big[\ln|\ln(-t)|\, \big]
 ^{\frac 34} \asA t \to 0^-,
  }
  \ee
  where $C>0$ is arbitrary.
Of course, unlike the
linear case in \ef{16}, the constant $C$ in \ef{m5}
 cannot play a role
 in view of the scaling invariance of the PDE \ef{u1},
  \be
  \label{sc1}
  u = A \hat u, \quad x =\sqrt{ A}\,  \hat x \quad (A>0),
  \ee
  which changes the constant $C$ without changing the regularity
  of the point. Similar to \ef{kk2}--\ef{kk4}, the scaling
  \ef{sc1} itself
 implies \ef{m5} by choosing $A(t)=a_0(\t)$ from \ef{m4}.


 It is not  difficult to perform a more detailed study of the
 original ODE \ef{m2}, though
 a too thorough   analysis seems  excessive.
  As usual, for the regularity,  an oscillatory cut-off
  is necessary. Thus, it turns out that for the quasilinear \ef{u1} and linear
\ef{ir12} fourth-order parabolic  equations, the regularity
conditions (including the concepts and the methods) of boundary
regularity analysis can be obtained in similar lines.

 An analogous study  can be done for the
 non-divergent {\em thin film equation} (TFE--4)
  $$
  u_t=-(u^2 u_{xxx})_x,
   $$
 though its existence-uniqueness theory is less developed, so the
 results will be more formal. Of course, the same techniques apply
 to the TFE--(2,2), with other distributions of the inner and
 outer derivatives in the differential form
  $$
  u_t= -(u^2 u_{xx})_{xx} ,
 $$
 and to other PDEs with more general nonlinearities. Extensions to
 various quasilinear counterparts of the $2m$th-order parabolic equations
 \ef{PE1} are also possible on the basis of the non self-adjoint spectral theory
 in \cite{Eg4}.

 \subsection{Fourth-order hyperbolic equation: regularity via Hermitian spectral theory for a
 pencil}
 \label{S4FF}

We recall that operator pencil theory fully occurred already in
the regularity study for elliptic PDEs  \cite{Kond67}, and was
later developed and extended in many papers;
see \cite{KondOl83, MazPl81, KMR1, KozMaz01, Freh06} for
references. Boundary regularity for hyperbolic PDEs is definitely
less developed; see last pages of Kondrat'ev--Oleinik's survey
\cite{KondOl83}  of 1983  for some references. Note that, for
fixed standard characteristic paraboloids posed for \ef{l1zz}, the
regularity problem falls into the scope of the results of
Kondrat'ev's ``parabolic" paper \cite{Kond66}. The case of
``expanding paraboloids" as $t \to 0^-$,  similar to the above
cases for the bi-harmonic equations and others, was not treated
before by the same reasons (the lack of proper spectral theory and
matching via a boundary layer approach, which are very difficult
to fully justify).

Here,  according to
the principle
 \ef{OscInf}, we very briefly discuss a new
types of operator pencils that is needed to tackle other
regularity problems. The simplest such model satisfying
\ef{OscInf} is the {\em linear fourth-order hyperbolic} ({\em
wave}) {\em equation}
\be
\label{l1zz}
 u_{tt}=-u_{xxxx} \inB Q_0,
 \ee
 also  known as the 1D {\em linear beam equation}.
We consider the same Dirichlet problem as in \ef{ir12} with two
bounded initial functions $u(x,-1)=u_0(x)$ and $u_t(x,-1)=u_1(x)$.
Concerning the fundamental solution  and  necessary spectral
properties, we follow \cite{GalpSWE}.

Thus, the {\em fundamental solution}
of (\ref{l1zz}) has
the  self-similar form
 \be
 \label{l2zz}
  \mbox{$
 b_0(x,t)= \sqrt t \,\,
  F(y), \quad y= \frac x{\sqrt t} \whereA b_0(x,0)=0, \quad b_{0t}(x,0)=\d(x)
   $}
  \ee
in the sense of bounded measures.
  The rescaled kernel $F_0=F_0(|y|)$ is symmetric and solves the ODE
   \be
   \label{l3zz}
   \textstyle{
  {\bf B}F \equiv  -F^{(4)} - \frac 14 \, F'' y^2 - \frac 14 \, F'y + \frac 14\, F=0
   \quad \mbox{in} \quad \re, \quad \int F=1.
 }
 \ee
 Integrating (\ref{l3zz}) once yields
 \be
   \label{l4}
   \textstyle{
  - F''' - \frac 14 \, F' y^2  + \frac 14\, F y=0 \LongA
 F(y)= \frac 1{2 \pi}\, \int\limits_0^\infty \frac
{\sin z \cos( \sqrt z y)}{z^{3/2}}\, {\mathrm d}z.
 }
  \ee
A WKBJ-type asymptotic analysis of the ODE (\ref{l4}) yields the
behaviour like  \ef{OscInf}:
 \be
 \label{l7}
  \textstyle{
 F(y) \sim C_1 y^{-\frac{11}{13}}\cos \big( \frac {y^2}4+C_2\big) \quad \mbox{as}
 \quad y \to +\infty, \quad C_1 \not = 0.
 }
  \ee

Similar to \ef{ph3}, the first scaling is
 \be
 \label{p1*}
 \begin{matrix}
 u(x,t)= v(y,\t), \quad y=\frac x{\sqrt {-t}}, \quad \t=-\ln(-t),
 \qquad\,\,\,
\qquad \ssk\ssk\ssk\\
 v_{\t\t}+ v_\t + v_{\t y}y= {\bf B}^* v  \equiv
 -v_{yyyy} - \frac 14 \, y^2 v_{yy} -
\frac{3}4\, y v_y.\qquad
  \end{matrix}
 \ee
Not paying enough attention to the functional setting of the
operators involved in weighted $L^2$-spaces  \cite{GalpSWE}, we
concentrate on the polynomial eigenfunctions of a quadratic pencil
to appear;
 see Markus \cite{Markus} for necessary
 concepts and theory of linear operator pencils. Namely, to find
eigenfunctions, we set
 \be
 \label{p10}
  w(y,\t)={\mathrm e}^{\l_k \t} \psi_k^*(y) \LongA
  {\bf C}^*(\l_k)\psi_k^* \equiv  {\bf B}^*\psi_k^*-(\l_k^2 + \l_k) \psi_k^*
  - \l_k (\psi^*_k)' \, y=0.
 \ee
Looking for finite polynomial eigenfunctions (recall that one
needs two sets of eigenfunctions $\Psi^*=\{\psi_k^*(y), \,
\phi_k^*(y)\}$ as well as the adjoint ones $\Psi=\{\psi_k(y), \,
\phi_k(y)\}$ \cite{GalpSWE}),
  \be
    \label{p110}
     \tex{
    \psi_k^*(y), \,\,\, \phi_k^*(y) =
    y^k+... \quad \mbox{are $k$th-order
    polynomials,}
     }
     \ee
  substituting into (\ref{p10}) and keeping the higher-degree
 terms
 $\sim y^k$ yields the following quadratic equation for
 eigenvalues:
 \be
 \label{llk}
  \mbox{$
  \l_k^2 + (k +1) \l_k + \frac{k(k-1)}4 + \frac{3k}4
  =0.
   $}
   \ee
   This gives two, shifted by -1, series of eigenvalues,
    \be
    \label{p11}
     \mbox{$
     \l_k^+=- \frac k2 \andA \l_k^-=-\frac k2 -1 \forA
     k=0,1,2, ... \, .
      $}
      \ee
Let us  specify a {\em half} of those {\em generalized Hermite
polynomials}, which is a complete closed set. There is another one
obtained in a similar manner, which is obviously necessary for the
second-order in $t$ PDE \ef{l1zz}.

\begin{proposition}
 \label{Pr.Eig*}
The  eigenfunctions of the pencil $(\ref{p10})$
 with the spectrum $\s_+=\{\l_k^+\}$ in  $(\ref{p11})$ are given
 by the normalized polynomials
 \be
  \label{p20NN}
   \mbox{$
   \psi_k^*(y)= \frac 1{\sqrt{ k!}} \, \big[y^k
   + \sum_{j=1}^{[ \frac k4]} \frac 1{3^j \, j!} \,
   (y^k)^{(4j)} \big], \quad k=0,1,2,... \, .
 $}
  \ee
 \end{proposition}

Assuming as usual that $\var(\t)$ is a slow growing function (so
\ef{vv1}, \ef{vv2},\,... hold), we introduce the variable $z$ in
\ef{z1} and the resulting BL-variables as in \ef{z2}, where
 \be
 \label{z2zz}
 \tex{
 \xi=\var^2(\t)(1-z), \quad z= \frac y{\var(\t)}, \quad D_y=-\var(\t) D_\xi.
 }
  \ee
  Assuming that the class of generic solutions admitting a
  BL-representation is under scrutiny and performing the cut-off
  as in
  \ef{a1}, where $\hat v_{\t\t}=v_{\t\t} H$, by the expansion \ef{a7} (we currently
   take into account the first coefficient only, with $\langle 1,F \rangle=1$ by
  the orthonormality, over the first half \ef{p20NN} of all the
   polynomials, and this defines the class of solutions involved), so  eventually, by the asymptotics \ef{l7} yields
   \be
   \label{a011}
 \tex{
   a_0''+ a_0'+... \sim  \g_1 C_3 a_0 \var^{\frac{28}{13}}\, \cos\big(\frac
   {\var^2}4+C_4\big)+...\quad \big(\frac{28}{13}=3- \frac
   {11}{13}\big),
 }
    \ee
 where we omit  all smaller terms containing $\var', \, \var''$, and
others. For slow varying coefficients $a_0(\t)$, the leading term
is indeed $a_0'$, so, in this critical case, the governing ODE is
\be
   \label{a011z}
 \tex{
   \frac { a_0'}{a_0} \sim \hat \g\, C_3  \var^{\frac{28}{13}}\, \cos\big(\frac
   {\var^2}4+C_4\big)+... \LongA \ln|a_0(\t)| \sim \int\limits^\t \var^{\frac{28}{13}}\, \cos\big(\frac
   {\var^2}4+C_4\big)\, {\mathrm d}\t.
 }
    \ee
Clearly, unlike \ef{14N}, by setting $\var^2(\t)=s$, we see from
\ef{vv1}  that a converging integral for the irregularity is
impossible. However, this does not mean a ``total" regularity.
 Actually, this again implies that the regularity issue here is rather
 subtle and, as usual for oscillatory kernels, can be achieved
 by a proper oscillatory cut-off of  suitable slow growing
 $\var(\t)$.
 For general
$\var(\t)$, a more delicate analysis of the ODE \ef{a011} with all
the terms included is necessary, which leads to some technical
questions.
 A
full justification of the general approach then becomes difficult
and even questionable.

\smallskip

{\bf Acknowledgement.} The author would like to thank M.~Kelbert
for a fruitful  discussion of the regularity-probability PDE
issues and interesting historical comments on Khinchin's results
in 1924 and  1936 \cite{Khin36}. The author also thanks
V.G.~Maz'ya for help concerning general concepts of boundary
regularity theory and for pointing out Fe${\rm \check{i}}$gin's
earlier results \cite{Fei71}, S.R.~Svirshchevskii for getting
early and rare V.P.~Miha${\rm \check{i}}$lov's papers (in Russian)
in 1961--63, and finally to A.~Kyprianou for a short but
informative conversation on Petrovskii's and Khinchin's
contributions to the LIL and other probability issues.


\begin{appendix}
\section*{Appendix A. Comment: On using majorizing order-preserving operators}
 \label{S3M}
 \setcounter{section}{1}
\setcounter{equation}{0}

\begin{small}


It is clear that the main difficulty in establishing optimal
regularity/irregularity criteria for higher-order linear parabolic
PDEs such as \ef{ir12}, \ef{PE1} (or third-order dispersion ones
in \ef{l1} and \ef{sss1} to be studied later on) is the fact that
these flows do not obey the Maximum Principle and do not exhibit
order-preserving features. As is well-known, the latter ones have
been key for
 the heat equation \ef{ir1} and other nonlinear second-order
 diffusion-reaction-absorption-convection-... equations; see
 references in \cite{Abd00, Abd05, Herr04, Moss00}.

Therefore,
 before going into eigenfunction expansion techniques for solving
 the regularity problem for \ef{ir12}, we cannot avoid a
 temptation to
  briefly discussing
 a rather formal opportunity to deal
  with the traditional
 techniques based on the Maximum Principle even for the
 higher-order parabolic flow.
 This is about the idea  of {\em majorizing
 order-preserving flows} for higher-order parabolic equations
 \cite{GPMaj}, which we explain relative to the poly-harmonic flow \ef{PE1}.
  Namely, given the rescaled oscillatory kernel \ef{FundSol}, we
 construct the corresponding majorizing kernel $\bar F(y)>0$ such that
  \be
  \label{FF1}
   \tex{
  |F(y)| \le D_m \bar F(y) \inB \ren \whereA \int \bar F=1,
 }
   \ee
and $D_m>1$  is a constant called the {\em order deficiency} (or
{\em order defect}) of this majorizing kernel. It is clear that
the optimal order deficiency is ($D_1^*=1$ since $F>0$)
 \be
 \label{FF2}
  \tex{
  D_m^*= \int |F|>1, \quad \mbox{corresponding to the non-smooth kernel}
   \quad \bar F(y) = \frac 1{D_m^*} \, |F(y)|
    }
     \ee
(sometimes, dealing with not $C^\iy$ kernels $\bar F(y)$ is not
convenient, especially when determining its
eigenfunctions by differentiation via \ef{eigen}, with $F \mapsto
\bar F$, etc.).

We then introduce the corresponding {\em majorizing
order-preserving integral evolution equation} in $\ren \times
\re_+$ (it does not have a PDE representation in general) defined
by the formal convolution
 \be
 \label{F3}
  \tex{
 \bar u(t) = \bar b(t) * \bar u_0, \,\,\mbox{with}\,\, \bar u_0(x) \ge 0 \whereA \bar b(x,t)= t^{-\frac 1{2m}}
 \bar F(y), \quad y= \frac x{t^{1/2m}}.
 }
  \ee

The corresponding {\em comparison theorem} is as follows
\cite{GPMaj}:
 \be
 \label{F4}
  \tex{
  |u_0(x)| \le \frac 1{D_m}\,  \bar u_0(x) \inB \ren \LongA |u(x,t)| \le \bar
  u(x,t) \inB \ren \times \re_+.
  }
   \ee
The proof is straightforward in view of \ef{FF1} via comparison of
the solution of \ef{F3} with the parabolic flow
 \be
 \label{F5}
  \tex{
 u(x,t) =  b(x-\cdot,t) *  u_0(\cdot) \equiv  t^{-\frac 1{2m}}\,\int F\big(
 \frac {x-y}{t^{1/2m}}\big)\, u_0(y)\, {\mathrm d}y
  \inB \ren \times \re_+.
  }
  \ee
 Indeed, subtracting \ef{F3} and \ef{F5} yields for the difference
 the desired inequality:
 \be
 \label{F5N}
  \begin{matrix}
w(x,t)= \bar u(x,t)-u(x,t)= t^{-\frac 1{2m}} \int \big[ \bar
F(\cdot)\, \bar u_0 - F(\cdot)\, u_0\big] \ssk\ssk\ssk \\
 \ge  t^{-\frac
1{2m}} \int \big[ \bar F(\cdot)\, \bar u_0 - \frac 1{D_m}\,
|F(\cdot)|\, \bar u_0\big]
   =t^{-\frac 1{2m}} \int \big[  \frac {D_m  \bar
F(\cdot) - |F(\cdot)|}{D_m} \, \bar u_0\big]  \ge 0,
 \end{matrix}
 \ee
 by the definition of the majorizing kernel \ef{FF1}.
Since $\bar F>0$, i.e., the integral majorizing evolution \ef{F3}
is order-preserving,  solutions $\bar u(x,t)$ admit constructing
super-solutions $\bar U(x,t)$ in the standard sense, i.e., the
following comparison holds:
 \be
 \label{F6}
  \bar U(t) \ge \bar b(t) * \bar u_0
  \LongA 0 < \bar u(x,t) \le \bar U(x,t) \inB \ren \times \re_+.
   \ee
In \cite{GPMaj}, such comparison ideas resolved some asymptotic
problems for  semilinear parabolic equations related to Fujita
critical exponents and around.

It is clear that to construct a majorizing comparison theory for
the IBVPs in domains such as $Q_0$, further estimates of ``Green's
function"  corresponding  to the ``majorizing non-translational
flow"
 \ef{F3} in $Q_0$
(as usual, such an estimate is expected to be related to \ef{FF1})
are necessary. These are difficult questions, which are not
discussed here.
We just mention that
 such a true possibility would allow, in a reasonable simpler (and rigorous!) way,
 to get a proper {\em upper} estimates on
   solutions in $Q_0$ to prove that $u(0,0)=0$, which would guarantee its {\em
   regularity}, by constructing  proper super-solutions. As a formal clue, it would be very efficiently and attractively easy
    to use Petrovskii's ansatz \ef{sub1},
    \be
    \label{jjj1}
   \tex{
   \bar U(x,t) \sim - \frac 1{|\ln(-t)|^{1+\e_1}}\, {\mathrm
   e}^{{-d_0 |x|^{\frac{2m}{2m-1}}}/{(-t)^{\frac 1{2m-1}}}} + \frac 1{\ln|\ln(-t)|},
    }
 \ee
 where we have replaced the Gaussian structure in \ef{sub1} by
that in the estimate \ef{es11}, which is supposed to be true for
the majorizing positive kernel $\bar F(y)$.
 As we will show (cf. \ef{16} for $m=2$),  the space-time structure \ef{jjj1} would  give
  a correct
 critical behaviour for $R(t)$ as $t \to 0^-$, but we do not know if this can
  be justified in such a comparison way.
  (an open problem, but an extremely attractive one).
  We must admit that,
 besides some other fundamental open questions, in any case,
 unlike \ef{sub1}, using \ef{jjj1} for any $m \ge 2$ assumes
  complicated {\em integral calculus} to establish the
 validity of the integral inequality in \ef{F6}.
 More importantly, such an integral majorizing
 approach cannot be used for any lower bound on $u(x,t)$ to get
 the irregularity counterpart, when $|u(0,0)|>0$. As one can see,
 corresponding alternative ``minorizing" theory in similar lines cannot be
 constructed in principle. Therefore, we had no intention and in
 fact
 necessity to check whether functions such as \ef{jjj1} can
 satisfy \ef{F6}  or similar integral inequalities.
 However, we state the following principal
 {\em open problem}\footnote{Personally, the author
 suspects that this is not possible, especially since the ``cut-off" must be then
   inherited by the method (how
 and why?);
 though then proving the negative result (a principal non-applicability of such barrier techniques
 for $m \ge 2$)
 would be also important.}:
 {\em Is there any hope to get sharp regularity results by using
  majorizing-like
   calculus associated with some integral operator
  in $Q_0$ with a positive kernel?}


\end{small}

\end{appendix}

\end{document}